\documentclass
{article}
\usepackage{amsmath}
\usepackage{amsfonts}
\usepackage{amssymb}
\usepackage{xspace}
\usepackage[title]{appendix}
\usepackage{xcolor}
\usepackage{graphicx}
\usepackage{graphics,color}
\usepackage{amsmath}
\usepackage{amssymb}
\usepackage{mathrsfs}
\usepackage{cite}
\usepackage{verbatim}
\usepackage{float}
\usepackage{amsthm}
\usepackage{textcomp}
\usepackage{subfig}
\usepackage{enumerate}
\usepackage{hyperref}

\date{}
\makeatletter \@addtoreset{equation}{section} \makeatother

\newenvironment{proofof}{\par\noindent{\sc Proof}
}{\hfill\llap{$\Box$}\vspace{1\baselineskip}\par\noindent}
\newtheorem{theorem}{Theorem}[section]
\newtheorem{proposition}[theorem]{Proposition}
\newtheorem{lemma}[theorem]{Lemma}
\newtheorem{corollary}[theorem]{Corollary}
\newtheorem{remark}[theorem]{Remark}
\newtheorem{definition}[theorem]{Definition}
\newtheorem{example}[theorem]{Example}
\newcommand{\beq}{\begin{equation}}
\newcommand{\eeq}{\end{equation}}

\newcommand{\ba}{\begin{array}}
\newcommand{\ea}{\end{array}}
\newcommand{\bt}{\begin{theorem}}
\newcommand{\et}{\end{theorem}}
\newcommand{\bp}{\begin{proposition}}
\newcommand{\ep}{\end{proposition}}
\newcommand{\bl}{\begin{lemma}}
\newcommand{\el}{\end{lemma}}
\newcommand{\bc}{\begin{corollary}}
\newcommand{\ec}{\end{corollary}}
\newcommand{\bi}{\begin{itemize}}
\newcommand{\ei}{\end{itemize}}
\newcommand{\ben}{\begin{enumerate}}
\newcommand{\een}{\end{enumerate}}
\newcommand{\bpf}{\begin{proof}}
\newcommand{\epf}{\end{proof}}
\newcommand{\bpff}{\begin{proofof}}
\newcommand{\epff}{\end{proofof}}
\newcommand{\bdf}{\begin{definition}\rm}
\newcommand{\edf}{\end{definition}}
\newcommand{\br}{\begin{remark}\rm}
\newcommand{\er}{\end{remark}}
\newcommand{\bex}{\begin{example}\rm}
\newcommand{\eex}{\end{example}}
\def\qu#1{{[\![\, #1\, ]\!]}}
\def\pri{\hbox to 10pt{\hfil\hbox to 0.4pt{\vrule height5pt width0.4pt
                 depth0pt}\vrule width5pt height0.4pt depth0pt\hfil}}
\newcommand{\gM}{{\bf M}}

\newcommand{\R}{{\cal R}}

\newcommand{\calL}{{\cal L}}

\newcommand{\gR}{{\mathbb R}}

\newcommand{\gZ}{{\mathbb Z}}

\newcommand{\Nat}{{\mathbb N}}

\newcommand{\A}{{\cal A}}
\newcommand{\D}{{\cal D}}

\newcommand{\G}{{\cal G}}

\newcommand{\M}{{\cal M}}
\newcommand{\Ha}{{\cal H}}

\newcommand{\cart}{\mathop{\rm cart}\nolimits}
\newcommand{\BV}{BV}

\newcommand{\rP}{{\mathbb P}}

\def\Lip{\mathop{\rm Lip}\nolimits}

\def\spt{\mathop{\rm spt}\nolimits}

\def\dist{\mathop{\rm dist}\nolimits}

\def\Int{\mathop{\rm int}\nolimits}

\def\det{\mathop{\rm det}\nolimits}

\def\sing{\mathop{\rm sing}\nolimits}

\def\sing{\mathop{\rm sing}\nolimits}
\newcommand{\If}{\ \mbox{\rm if}\ }

\newcommand{\Id}{\mbox{\rm Id}}
\newcommand{\wc}{\rightharpoonup}

\def\Lip{\mathop{\rm Lip}\nolimits}

\newcommand{\Sph}{{\mathbb S}}

\newcommand{\gm}{{\bf{m}}}

\def\Det{\mathop{\rm Det}\nolimits}
\def\adj{\mathop{\rm adj}\nolimits}

\def\DDiv{\mathop{\rm Div}\nolimits}
\def\BVc{\stackrel{\BV}{\rightarrow}}

\let\a=\alpha
\let\be=\beta

\let\d=\delta
\let\e=\varepsilon

\let\vf=\varphi

\let\l=\lambda
\let\m=\mu
\let\n=\nu
\let\o=\omega
\let\p=\pi
\let\q=\psi
\let\r=\rho
\let\s=\sigma
\let\t=\theta
\let\tt=\tau
\let\x=\xi
\let\y=\eta

\let\vf=\varphi
\let\DD=\Delta

\let\SS=\Sigma

\let\wih=\widehat
\let\wdg=\wedge

\let\wid=\widetilde
\let\join=\bowtie
\let\pa=\partial
\let\sb=\subset

\let\emp=\emptyset

\let\fa=\forall
\let\tim=\times
\let\mx=\mbox
\let\bcup=\bigcup

\let\sm=\setminus
\let\ol=\overline

\let\ds=\displaystyle

\let\lan=\langle
\let\ran=\rangle

\let\i=\infty
\let\lm=\limits
\let\nolm=\nolimits
 \title{Strict BV relaxed area of Sobolev maps into the circle:\\
the high dimension case}
\date{\today}

\author{
Simone Carano\footnote{
Area of Mathematical Analysis, Modelling, and Applications,
             Scuola Internazionale Superiore di Studi Avanzati ``SISSA'',
Via Bonomea, 265 - 34136 Trieste, Italy. E-mail: scarano@sissa.it
                         }
                          \and
Domenico Mucci\footnote{Dipartimento di Scienze Matematiche, Fisiche e Informatiche, Universit\`{a}
di Parma, Parco Area delle Scienze 53/A - 43124
Parma, Italy.
E-mail: domenico.mucci@unipr.it
}
}

\begin{document}
\maketitle
\begin{abstract}
We deal with the relaxed area functional in the strict $BV$-convergence of non-smooth maps defined in domains of generic dimension and taking values into the unit circle. In case of Sobolev maps, a complete explicit formula is obtained. Our proof is based on tools from Geometric Measure Theory and Cartesian currents. We then discuss the possible extension to the wider class of maps with bounded variation. Finally, we show a counterexample to the locality property in case of both dimension and codimension larger than two.
\end{abstract}

\vspace{2mm}

\noindent {\bf AMS (MOS) 2022 subject clas\-si\-fi\-ca\-tion:} 49J45, 49Q05, 49Q20, 28A75.

%
\section*{Introduction}
In this paper, we continue the analysis of the explicit formula for the {\em relaxed area} functional with respect to the {\em strict convergence} for non-smooth vector-valued functions $u:B^n\to\gR^2$.

For a smooth function $u:B^n\to\gR^2$, we denote by $\mathcal{A}(u,B^n)$ the area of the graph of $u$, given by
$$ \mathcal{A}(u,B^n):=\int_{B^n}\sqrt{1+|\nabla u|^2+|M_2(\nabla u)|^2}\,dx  $$
where $|M_2(\nabla u)|^2$ is the sum of the square of all $2\tim 2$ minors of the gradient matrix $\nabla u$, so that $|M_2(\nabla u)|=|\det\nabla u|$ if $n=2$. In the sequel we shall write simply $\mathcal{A}(u)=\mathcal{A}(u,B^n)$.

Working with the natural $L^1$-convergence, the relaxed area functional is defined for every summable function $u\in L^1(B^n,\gR^2)$ by
\begin{align}\label{L1_area}
 \ol{\mathcal{A}}_{L^1}(u):=\inf\Bigl\{\liminf_{k\to\i}\mathcal{A}(u_k)\mid \{u_k\}\sb C^1(B^n,\gR^2),\,\,u_k\to u\,\,{\text{strongly in}}\,\,L^1\Bigr\}\,.
\end{align}
If $\ol{\mathcal{A}}_{L^1}(u)<\i$, then necessarily $u$ is a function of {\em bounded variation}. However, even in low dimension $n=2$, it turns out that the localized functional $A\mapsto \ol{\mathcal{A}}_{L^1}(u,A)$ fails to be subadditive, and hence it cannot be extended to a Borel measure. This behavior was conjectured by De Giorgi in \cite{DG} and proved by Acerbi and Dal Maso in \cite{AcDM}, where it is shown that non-subadditivity phenomena arise even for very simple cases like the vortex map $u_V(x)=x/|x|$ and the symmetric triple junction map $u_T$. A precise computation of the values $\ol{\mathcal{A}}_{L^1}(u_V)$ and $\ol{\mathcal{A}}_{L^1}(u_T)$ can be found in \cite{BES}(see also \cite{BMS}) and \cite{BP,S} respectively.  Moreover, for the analysis of the triple junction map without symmetry assumptions, we refer to \cite{BEPS}, where the authors provide an upper bound for the respective $L^1$-relaxed area \eqref{L1_area}, conjectured to be optimal. Other interesting upper bounds were recently obtained in \cite{BeScSc} for Sobolev maps valued in $\Sph^1$ and in \cite{ScSc} for piecewise constant maps taking three values.

\smallskip
The non-locality feature previously outlined makes quite challenging the relaxation analysis of $\mathcal A$. For this reason, it is interesting to
consider some variants of  \eqref{L1_area}, for example
by strengthening the topology of the convergence of $u_k$ to $u$ (see \cite{BePaTe:15,BePaTe:16,DeP,GMSl2}).
In recent years it has been proposed to impose the strict $BV$-convergence.
Referring to e.g. \cite{AFP} for the notation adopted in this paper, we only recall here that a sequence $\{u_k\}\sb \BV(B^n,\gR^2)$ is said to
converge to $u\in \BV(B^n,\gR^2)$ strictly in the $\BV$-sense, say $u_k\BVc u$, if $u_k\to u$ in $L^1(B^n,\gR^2)$ and $|Du_k|(B^n)\to|Du|(B^n)$, where $|Du|$ denoted the total variation of the distributional derivative $Du$.

For $u\in\BV(B^n,\gR^2)$, we thus denote
\beq\label{Area-rel} \ol{\mathcal{A}}_{\BV}(u):=\inf\Bigl\{\liminf_{k\to\i}\mathcal{A}(u_k)\mid \{u_k\}\sb C^1(B^n,\gR^2),\,u_k\BVc u\,{\text{strictly in}}\,\BV\Bigr\}.
\eeq
The reason for this choice is that we expect that whenever $\ol{\mathcal{A}}_{\BV}(u)<\i$, then the localized functional $A\mapsto \ol{\mathcal{A}}_{\BV}(u,A)$ gives rise to a Borel measure.
In case of low dimension $n=2$, partial results concerning the explicit formula of $\ol{\mathcal{A}}_{\BV}(u)$ have been obtained in \cite{BSS,BSS2,Ca}.

In this paper, we focus on the class of $\BV$-maps taking values into the unit circle.
More precisely, we denote by
$$\Sph^1:=\{y\in\gR^2\,:\,|y|=1\}\,,\quad D^2:=\{y\in\gR^2\,:\,|y|<1\} $$
the unit circle and disk in the target space, and for $X=\BV$ or $W^{1,1}$, we let
$$ X(B^n,\Sph^1):=\{u\in X(B^n,\gR^2)\mid |u(x)|=1\quad{\mx{\text for $\calL^n$-a.e. }}x\in B^n\} $$
where $\calL^n$ is the Lebesgue measure, and we focus on the class of Sobolev maps $W^{1,1}(B^n,\Sph^1)$.

In low dimension $n=2$, the following result was obtained in \cite{BSS}.
\bt Let $u\in W^{1,1}(B^2,\Sph^1)$, and let $\Det\nabla u$ denote the {\em distributional determinant} of $u$. Then,
$$ \ol{\mathcal{A}}_{\BV}(u)<\i \iff |\Det\nabla u|(B^2)<\i\,. $$
In that case, moreover, one has:
$$ \ol{\mathcal{A}}_{\BV}(u) =\int_{B^2}\sqrt{1+|\nabla u|^2}\,dx+|\Det\nabla u|(B^2)\,. $$
\et
\par The previous result says that the {\em energy gap} is detected by the distributional determinant.
Referring to the next section for the notation adopted here, we recall that in any dimension $n\geq 2$, the {\em current carried by the graph} of a Sobolev map
$u\in W^{1,1}(B^n,\Sph^1)$ is an integer multiplicity rectifiable $n$-current $G_u$ in $B^n\tim\Sph^1$, with finite mass
$$ \gM(G_u)=\int_{B^n}\sqrt{1+|\nabla u|^2}\,dx<\i\,. $$
Moreover, the relevant singularities of $u$ are detected by the current $G_u$, i.e., they can be described by homological arguments.

More precisely, denoting by $\o_2$ the closed 1-form in $\Sph^1$ 
$$ 
\o_2:=\ds \frac 12\,\bigl(y^1\,dy^2-y^2\,dy^1 \bigr) $$
the singularities are read by the $(n-2)$-dimensional current $\rP(u)\in\D_{n-2}(B^n)$ defined by
$$ \rP(u)(\y):= -\frac 1\pi G_u(d\y\wedge\o_2)\,,\quad \y\in \D^{n-2}(B^n)\,. $$
Therefore, in low dimension $n=2$ one has
$$ \pi\cdot\rP(u)(\y)=\langle\Det\nabla u,\y\rangle\quad\fa\,\y\in C^\i_c(B^2)\,. $$
If e.g. $u(x)=x/|x|$, one gets $\rP(u)=\d_{0_{\gR^2}}$, the unit Dirac mass at the origin.

\smallskip\par

In our Main Result, we extend the previous explicit formula to any high dimension $n$.
\bt\label{MainResult}
Let $n\geq 2$ and $u\in W^{1,1}(B^n,\Sph^1)$. Then, $\ol{\mathcal{A}}_{\BV}(u)<\i$ if and only if the $(n-2)$-current $\rP(u)$ is i.m. rectifiable and with finite mass, $\gM(\rP(u))<\i$.
In that case, moreover, one has:
$$ \ol{\mathcal{A}}_{\BV}(u) =\int_{B^n}\sqrt{1+|\nabla u|^2}\,dx+\pi\,\gM(\rP(u))\,. $$
\et
\par For our purposes, we exploit in our context a result taken from \cite{Mu22}. It says that if $u\in W^{1,1}(B^n,\Sph^1)$ satisfies $\ol{\mathcal{A}}_{\BV}(u)<\i$, then there exists a unique optimal {\em Cartesian current} $T_u$ that encloses the graph of $u$, and it is given by
$$ T_u=G_u+(-1)^{n-2}\rP(u)\times\qu{D^2}\,. $$
Therefore, the proof of the energy lower bound readily follows: by Federer's closure-compactness theorem, for every smooth sequence $\{u_h\}$ strictly converging to $u$, the graph $G_{u_h}$ weakly converges to $T_u$ (up to extracting a subsequence), and one concludes from the semicontinuity of the mass.
On the other hand, the energy upper bound holds true as a consequence of the following approximation result:
\bt\label{Tappr-int}
Let $n\geq 2$ and $u\in W^{1,1}(B^n,\Sph^1)$ be a Sobolev map with finite relaxed energy \eqref{Area-rel}. Then, there exists a smooth sequence
$\{u_h\}\sb C^\infty(B^n,\gR^2)$ such that $G_{u_h}\wc T_u$ weakly in $\D_n(B^n\times\gR^2)$ and
$\gM(G_{u_h})\to \gM(T_u)$ as $h\to\infty$.
\et
\par In Section~\ref{Sec:not}, we collect some notation and preliminary results.
In Section~\ref{Sec:example}, we give an explicit example, showing the strategy we follow in the proof of the relaxed formula.
In Section~\ref{Sec:main} we prove our Main Result, Theorem~\ref{MainResult}.
The proof of the approximation theorem~\ref{Tappr-int} is based on several technical results, the proof of which is postponed to Section~\ref{Sec:proofs}, for the sake of clarity. The fundamental step at the base of Theorem~\ref{Tappr-int} is contained in Theorem~\ref{TBet} and consists to reduce the proof to the case $u$ is smooth outside a ``nice'' singular set, precisely $\rP(u)$ is a polyhedral chain. This reduction can be done provided that the mass of the singularities current $\rP(u_k)$ of the modified map $u_k$ converges to the mass of $\rP(u)$. We point out that by a direct application of Bethuel's approximation theorem \cite[Thm.~2]{Bet} and Hardt-Pitts results in \cite{HP}, one obtains the flat norm convergence of $\rP(u_k)$ to $\rP(u)$, which is not enough for our purpose.
The actual proof requires a deeper use of Bethuel's result in a more involved construction argument, based on Federer's strong polyhedral approximation theorem. Once $u$ can be supposed to be smooth out of the support of the polyhedral $(n-2)$-chain $\rP(u)$, by a standard argument based on the dipole construction idea, we can build a recovery sequence for the energy (Theorem~\ref{Tdipolo}), taking care of removing higher codimension singularities generated in the dipole construction (Theorems~\ref{Tsing} and \ref{Tsimplex}).

Finally, in Section~\ref{Sec:final} we briefly discuss some related open questions, mainly concerning the validity of an explicit formula of the relaxed energy in the wider class of maps in $\BV(B^n,\Sph^1)$. Moreover, we show the non-locality of $\overline{\mathcal{A}}_{BV}$ in dimension and codimension greater than $2$. Precisely, the set function $A\to\overline{\mathcal{A}}_{BV}(u,A)$ fails to be subadditive for $u:B^3\to\gR^3$, even in the Sobolev case, as provided by the vortex map $u_V(x)=x/|x|$.
\section{Notation and preliminary results}\label{Sec:not}
In this section, we collect some background material and preliminary results.
\subsection{Functions of bounded variation}
Referring to \cite{AFP} for the notation on $\BV$-functions, we recall that $u$ belongs to $\BV(B^n,\gR^2)$ if $u\in L^1(B^n,\gR^2)$ and the distributional derivative $Du$ is an $\gR^{2\tim n}$-valued Borel measure with finite total variation. The usual decomposition
$$ Du=D^au+D^Cu+D^Ju $$
into the mutually singular absolutely continuous, Cantor, and Jump components is adopted. In particular, $D^au=\nabla u\,d\calL^n$, where $\nabla u\in L^1(B^n,\gR^{n\tim 2})$ is the approximate gradient and $\calL^n$ the Lebesgue measure.
The Jump component is given by $D^Ju=(u^+-u^-)\otimes\nu\,\Ha^{n-1}\pri J_u$, where $\Ha^{n-1}$ is the Hausdorff measure, $J_u$ is the Jump set of $u$, a countably $(n-1)$-rectifiable set of $B^n$, $\nu$ a unit normal to $J_u$, and $u^\pm$ the approximate limits of $u$ at points in $J_u$ w.r.t. the given unit normal.
The Cantor component satisfies $|D^Cu|(B)=0$ for each Borel set $B\sb B^n$ such that $\Ha^{n-1}(B)=0$.
Therefore, if $D^Ju=D^Cu=0$, actually $u$ is a Sobolev function in $W^{1,1}(B^n,\gR^2)$.
Finally, we recall that the {\em strict convergence} $u_k\BVc u$ in $\BV(B^n,\gR^2)$ is given by the strong convergence $u_k\to u$ in $L^1(B^n,\gR^2)$ joined
with the total variation convergence $|Du_k|(B^n)\to |Du|(B^n)$, as $k\to\infty$.
\subsection{Rectifiable currents}
For a given open set $U\sb\gR^N$, the space $\D_k(U)$ of $k$-dimensional currents in
$U$ is the strong dual of the space $\D^k(U)$ of compactly supported smooth $k$-forms in $U$, for $k=0,\ldots,N$.
For any $T\in\D_k(U)$, we define its \emph{mass}
$\gM(T)$ as
$$\gM(T):=\sup\{ T(\o) \mid \o\in\D^k(U)\,,\,\,\Vert\o\Vert\leq 1\}\,, $$
where $\Vert\o\Vert$ is the {\em comass norm}.\\
The {\em weak convergence} $T_h\wc T$ in $\D_k(U)$ is defined by the convergence
$$ \lim_{h\to\infty}T_h(\o)=T(\o)\quad\fa\,\o\in\D^{k}(U) $$
and in that case one has
$$ \gM(T)\leq\liminf_{h\to\infty}\gM(T_h)\,. $$
For $k\geq 1$, the {\em boundary} of $T\in\D_k(U)$ is the $(k-1)$-current $\partial T$ defined by relation
$$\pa T(\y) := T(d\y)\,, \qquad \y\in \D^{k-1}(U) $$
where ${d}\y$ is the differential of $\y$, and we set $\pa T=0$ if $k=0$.
For $k\geq 1$, a $k$-current $T$ with finite mass is called {\em rectifiable} if
there exist a $k$-rectifiable set $\M$ in $U$, an $\Ha^k\pri\M$-measurable function $\x:\M\to\Lambda^k\gR^m$ such that $\x(x)$ is a simple unit
$k$-vector orienting the approximate tangent space to $\M$ at $\Ha^k$-a.e. $x\in\M$, and an $\Ha^k\pri\M$-summable and non-negative function
$\t:\M\to[0,+\infty)$ such that
$$ T(\o)=\int_\M \t\,\langle \o,\x\rangle\,d\Ha^k \qquad\fa\,\o\in\D^k(U)\,.$$
We thus get
$\gM(T)=\int_\M \t\,d\Ha^k <\i$.
In addition, if the multiplicity function $\t$ is integer-valued, the current $T$ is called {\em i.m. rectifiable} and the corresponding class is denoted by $\R_k(U)$.
If e.g. $\M$ is a smooth $k$-manifold in $U$ with $\Ha^k(\M)<\i$, taking $\theta=1$ we obtain the current $\qu{\M}\in\R_k(U)$ whose action on $k$-forms agrees with the classical notation from Differential Geometry. In particular, for $k=0$, a current $T$ in $\R_0(U)$ is given by
$$T=\sum_{i=1}^m d_i\,\d_{a_i} $$
where $m\in\Nat^+$, $d_i\in\gZ$, $a_i\in U$ for $i=1,\ldots,m$, and $\d_{a}$ is the unit Dirac mass at a point $a\in U$.\\
Finally, a current $T$ is called \textit{integral} if both $T$ and $\partial T$ are i.m. rectifiable currents. By the boundary rectifiability theorem (cf. \cite[30.3]{Si}), if $T$ is i.m. rectifiable and $\gM(\partial T)<\infty$, then $T$ is integral.
We refer to \cite[Ch.~6]{Si} and \cite[Ch.~2]{GMSl1} for further details.
%
%
\subsection{Graph currents}
If $u\in C^1(\ol B^n,\gR^2)$, the graph current $G_u$ in $\R_n(B^n\tim\gR^2)$ is given by integration on the oriented graph $n$-manifold $\G_u$.
Therefore, by the area formula we equivalently have
\beq\label{Gu}
G_u(\o):=\int_{B^n}(\Id\join u)^\#\o\,,\qquad \o\in\D^n(B^n\tim \gR^2)
\eeq
where $(\Id\join u)(x):=(x,u(x))$ is the graph map, and its mass satisfies
\beq\label{MGu} \gM(G_u)=\Ha^n(\G_u)=\mathcal{A}(u)=\int_{B^n}\sqrt{1+|\nabla u|^2+|M_2(\nabla u)|}\,dx \,. \eeq
\par
To every Sobolev map $u\in BV(B^n,\Sph^1)$, we associate
the $n$-current $G_u$ in $\R_n(B^n\tim \gR^2)$ carried by the ``graph'' of $u$. It is given again by \eqref{Gu}, where this time the pull-back makes sense in terms of the approximate gradient $\nabla u$ of $u$.
Every $n$-form $\o\in\D^n(B^n\tim \gR^2)$ splits as $\o^{(0)}+\o^{(1)}+\o^{(2)}$ according to
the number of ``vertica'' differentials. Writing $\o^{(0)}=\phi(x,y)\,dx$ for some $\phi\in
C^\i_c(B^n\tim \gR^2)$, where $dx:=dx^1\wdg\cdots\wdg dx^n$, we have
$$ G_u(\phi(x,y)\,dx)=\int_{B^n}\phi(x,u(x))\,dx\,. $$
Setting moreover $\wih{dx^i}:=dx^1\wdg\cdots\wdg dx^{i-1}\wdg dx^{i+1}\wdg\cdots\wdg dx^n$, we may write
\beq\label{om}
\o^{(1)}=\sum_{i=1}^{n}\sum_{j=1}^{2}(-1)^{n-i}\phi^j_i(x,y)\,\wih{dx^i}\wdg dy^j
\eeq
for some $\phi^j_i\in C^\i_c(B^n\tim \gR^2)$, and we obtain
$$ G_u(\o^{(1)}) =  \ds\sum_{j=1}^2\sum_{i=1}^n
\ds\int\nolm_{B^n}\nabla_i
u^j(x)\phi^j_i(x,u(x))\,dx\,. $$
Finally, by the area formula we have
$$ G_u(\o^{(2)})=0\quad\fa\,\o\in\D^n(B^n\tim\gR^2)\,.
$$
In particular, we get
$$ \gM(G_u)=\int_{B^n}\sqrt{1+|\nabla u|^2}\,dx\,. $$
\par In general, the graph current
$G_u$ of a Sobolev map $u\in W^{1,1}(B^n,\Sph^1)$ has a non zero boundary in $B^n\tim \gR^2$.
Taking for example $n=2$ and $u(x)=x/|x|$, we have:
$$\pa G_u\pri B^2\tim \gR^2=-\d_{0_{\gR^2}}\tim\qu{\Sph^1}\,. $$
%
However, a density argument shows that the boundary current $\pa G_u$ is null on every $(n-1)$-form in $B^n\tim \gR^2$ which
has no ``vertical'' differentials. Moreover, $G_u$ is an integral flat chain in $B^n\tim\gR^2$ with support contained in $\ol B^n\tim\Sph^1$.
Therefore, by Federer's flatness theorem we can see $G_u$ as a current in $\R_n(B^n\tim\Sph^1)$, and actually
\beq\label{bdGu} (\pa G_u)\pri B^n\tim\gR^2=(\pa G_u)\pri B^n\tim\Sph^1\,.
\eeq
\subsection{Singularities}
If $u\in W^{1,1}(B^n,\gR^2)\cap L^\infty(B^n,\gR^2)$, it is well defined the distribution
\beq\label{Divmudef}
\DDiv_{\bar\a}\gm_u := \frac
1{2}\sum_{j=1}^2\sum_{i\in\ol\a}\frac{\pa}{\pa
x_i}\bigl(u^j(x)\,((\adj \nabla u)_{\ol\a})^j_i\bigr)  \eeq
for each ordered multi-index $\a$ of length $n-2$ in $\{1,\ldots,n\}$, where $\bar\a$ is the complementary ordered index of length two.  For $n=2$, the right hand side of definition \eqref{Divmudef} reduces to the distributional determinant $\Det\nabla u$.
In high dimension $n\geq 3$, instead, we obtain the $\bar\a$-component of the distributional Jacobian $J(u)$, which can be viewed as an $\gR^{d(n)}$-valued distribution, with $d(n)=n(n-1)/2$. The notion of distributional Jacobian was first introduced in \cite{Mo} (see also \cite{Re,Mul,Ba}) to analyse singularities of non-smooth maps and has been widely studied in the literature, together with the related notion of relaxed Jacobian total variation \cite{DL,FFM,Mu06,Mu10,Pa,DeP}. Notice that
$$
\DDiv_{\bar\a}\gm_u=M_2(\nabla u)_{\bar\a}\qquad\mbox{if }u\mbox{ is smooth},
$$
where $M_2(\nabla u)_{\bar\a}$ is the $2\tim 2$ minor of the gradient matrix $\nabla u\in\gR^{2\tim n}$ with columns detected by $\bar\a$.\\
The measure $\DDiv_{\bar\a}\gm_u$ can be defined also for any $BV$-map $u$ with finite relaxed energy \eqref{Area-rel}, by considering $Du$ in place of $\nabla u$, see \cite{Mu22} for further details.\\
Now suppose that $u\in W^{1,1}(B^n,\Sph^1)$. We can easily relate the distributional Jacobian $J(u)$ to an i.m. rectifiable current $\rP(u)$ defined as follows.
Let $\p:B^n\tim \gR^2\to B^n$ and $\wih\p:B^n\tim \gR^2\to
\gR^2$ denote the orthogonal projections onto the first and second
factor, respectively. The current $\rP(u)\in\D_{n-2}(B^n)$ of the singularities of $u$ is given by
\beq\label{Pu} \rP(u)(\y):= -\frac 1\pi G_u(\p^\#d\y\wedge\wih\p^\#\o_2)\,,\quad \y\in \D^{n-2}(B^n),
\eeq
where $\o_2$ denote the closed 1-form in $\Sph^1$
\beq\label{o2} \o_2:=\ds \frac 12\,\bigl(y^1\,dy^2-y^2\,dy^1 \bigr)\,, \eeq
so that $\o_2$ is a generator of the first cohomology group of $\Sph^1$, and
$d\o_2=dy:=dy^1\wdg dy^2$, as a form in $\gR^2$.
In the sequel, when it is clear from the context we omit to write the action of the projection maps.
Since $ d(\y\wedge\o_2)=d\y\wedge \o_2+(-1)^{n-2}\y\wedge dy$, whereas
$$ G_u(d\y\wedge \o_2)=-\pi\cdot\rP(u)(\y)\,,\quad G_u(\y\wedge dy)=0, $$
on account of \eqref{bdGu} we obtain:
\beq\label{bdGucont} (\pa G_u)\pri B^n\tim\gR^2= (\pa G_u)\pri B^n\tim\Sph^1=-\rP(u)\tim\qu{\Sph^1}\,.
\eeq
Moreover, by the very definition it turns out that
\beq\label{bdPu} \pa\rP(u)\pri B^n=0\,. \eeq
This property is trivial when $n=2$, whereas in high dimension $n\geq 3$, for every $\vf\in\D^{n-3}(B^n)$ we get
$$ \pa\rP(u)(\vf)=\rP(u)(d\vf)=-\frac 1\pi G_u(\p^\#d(d\vf)\wedge\wih\p^\#\o_2)=0, $$
since $d(d\vf)=0$.
\par For future use, we recall that a Cartesian current $T$ in $B^n\times \Sph^1$ with underlying map $u$ in $W^{1,1}(B^n,\Sph^1)$ is given by
\beq\label{Tcart} 
T=G_u+L\times\qu{\Sph^1} 
\eeq
for some i.m. rectifiable current $L\in\R_{n-1}(B^n)$, with finite mass, satisfying the boundary condition $(\pa L)\pri B^n=\rP(u)$, compare \cite[Sec.~1.5]{GM2} or 
\cite[Sec.~6.2.2]{GMSl2}.
\bex If $u\in W^{1,1}(B^2,\Sph^1)$, we have
$\pi\,\rP(u)=\Det\nabla u$. In particular, if $u$ is smooth outside a finite set of points $\SS=\{a_1,\ldots,a_m\}$, we obtain
\beq\label{Pu2} \rP(u)=\sum_{i=1}^m \deg(u,a_i)\,\d_{a_i},
\eeq
where $\deg(u,a_i)\in\gZ$ is the Brouwer degree\footnote{We refer to \cite[Sec.~3.1]{GMSl1} for the definition of Brouwer degree of a Sobolev map, see also \cite[Sec.~2.3]{BSS}. } of $u$ around the point $a_i$. For example, with $u(x)=x/|x|$, we get $\rP(u)=\d_{0_{\gR^2}}$.
In high dimension $n\geq 3$, for any $u\in W^{1,1}(B^n,\Sph^1)$ we get $\pi\,\rP(u)=J(u)$.
In Section~\ref{Sec:example}, we deal with the map
$$u(x)=\frac{\wid x}{|\wid x|}\,,\qquad x=(\wid x, \wih x)\in \gR^2\tim\gR^{n-2}\,,$$
so that $\rP(u)=(-1)^{n-2}\qu{\DD^{n-2}}$,
where $\qu{\DD^{n-2}}$ is the $(n-2)$-current given by integration on the naturally oriented $(n-2)$-disk
$$ \DD^{n-2}:=\{(0,0,\wih x)\in\gR^n\,:\, |\wih x|\leq 1\}\,. $$
\eex
\subsection{Stratification}
If $n\geq 3$, a current $T\in\R_n(B^n\tim\gR^2)$ is identified by the measures
$$ \ba{c} \m_h[T]:=T\pri dx\,,\qquad \m^j_i[T]:=T\pri (-1)^{i-1} dy^j\wdg \wih {dx^i}\,,
\\ %
 \m^{\bar\a}_v[T]:=T\pri \s(\a,\bar\a)\,dx^\a\wdg dy\,,\qquad dy:=dy^1\wdg dy^2 \ea $$
for each $i=1,\ldots,n$, $j=1,2$, and each ordered multi-index $\a$ of length $n-2$ in $\{1,\ldots,n\}$, where the sign $\s(\a,\bar\a)=\pm 1$ is such that $dx^\a\wdg dx^{\bar\a}=\s(\a,\bar\a)\,dx$.
We also fix an order on the set of the $d(n):=n(n-1)/2$ multi-indexes $\bar\a$ of length two in $\{1,\ldots,n\}$,
and we correspondingly denote by $\m_v[T]$ the $\gR^{d(n)}$-valued measure in $B^n\tim\gR^2$ with components $\m^{\bar\a}_v[T]$.
If $n=2$, then $\m_v[T]:=T\pri dy$.
\par Notice that if $T=G_u$ for some smooth function $u\in C^1(\ol B^n,\gR^2)$, by \eqref{Gu} we readily obtain $ \m_h[G_u]=(\Id\join u)_\#(\calL^n\pri B^n)$,
$\m^j_i[G_u]=(\Id\join u)_\#(\nabla_iu^j\,\calL^n\pri B^n)$, and also
$$ \m^{\bar\a}_v[G_u]=(\Id\join u)_\#(M_2(\nabla u)_{\bar\a}\,\calL^n\pri B^n) \qquad\fa\,\a.$$
\subsection{The optimal lifting Cartesian current}
Assume now that $u\in W^{1,1}(B^n,\Sph^1)$ has finite relaxed energy \eqref{Area-rel}.
Then, viewing $G_u$ as a current in $B^2\tim\gR^2$, by the results from \cite{Mu22}, it turns out that there exists a unique
i.m. rectifiable current $T_u\in\R_n(B^n\tim\gR^2)$ satisfying the following properties:
\ben
\item $\gM(T_u)<\infty$ and $(\pa T_u)\pri B^n\tim\gR^2=0$;
\item if $S_u:=T_u-G_u$, then $S_u$ is completely vertical, i.e., $S_u(\o)=0$ for every $\o\in\D^n(B^n\tim\gR^2)$ such that $\o^{(2)}=0$.
\een
In particular, $T_u$ is a Cartesian current in $\cart(B^n\tim\gR^2)$, see
\cite[Ch.~4]{GMSl1}, and
$$ \gM(T_u)=\gM(G_u)+\gM(S_u)\,. $$
\par
More precisely, the horizontal component of $T_u$ satisfying $\m_h[T_u]=(\Id\join u)_\#(\calL^n\pri B^n)$, we require that the intermediate components only depend on
$u$ through formulas
\beq\label{mijT} \m_i^j[T_u] =\m^j_i[u] \quad \fa\,i,j \eeq
where $\m^j_i[u]$ is the minimal lifting measure in the sense of Jerrard-Jung \cite{JJ}. Therefore,
$$ \m^j_i[u]=(\Id\join u)_\#(\nabla_iu^j\,\calL^n\pri B^n)\,. $$
For each multi-index $\a$ of length $n-2$ as above, we thus get
\beq\label{lift}  \ds\int_{B^n} g(x)\,d\m_v^{\bar\a}[T_u]=\lan\DDiv_{\bar\a}\gm_u,g\ran \qquad\fa\, g\in C^\i_c(B^n), \eeq
where $\DDiv_{\bar\a}\gm_u $ is defined in \eqref{Divmudef},
so that actually
\beq\label{PuDivmu}
\rP(u)(g(x)\,dx^{\a})=\frac 1\pi\,(-1)^{n-2}\s(\a,\bar\a)\,\lan \DDiv_{\bar\a}\gm_u ,g\ran\quad\fa\,g\in C^\i_c(B^n)\,. \eeq
\par We are now in position to prove the following
\bt\label{TPu} Let $n\geq 2$ and $u\in W^{1,1}(B^n,\Sph^1)$ be a Sobolev map with finite relaxed energy \eqref{Area-rel}. Then
\beq\label{Su} S_u=(-1)^{n-2}\rP(u)\tim\qu{D^2}
\eeq
where $\rP(u)$ is an i.m. rectifiable current in $\R_{n-2}(B^n)$ with finite mass and no inner boundary, see \eqref{bdPu}.
\et
\bpf By \eqref{lift}, for every $\a$ we get the total variation bound:
$$
 |\DDiv_{\bar\a}\gm_u|(B^n)\leq |\m_v^{\bar\a}[T_u]|(B^n\tim\gR^2)<\i \,. $$
As a consequence, equation \eqref{PuDivmu} implies that the current $\rP(u)$ has finite mass.

On the other hand, by \cite{GM} we already know that the relaxed total variation energy of $u$ as a map in $\BV(B^2,\Sph^1)$ is finite, whence the class of Cartesian currents in $B^n\times\Sph^1$ with underlying map equal to $u$ is non-empty, see \eqref{Tcart}.
Therefore, there exists $L\in\R_{n-1}(B^n)$ such that $(\pa L)\pri B^n=\rP(u)$, i.e., it turns out that $\rP(u)$ is an integral flat chain.
As a consequence, by the boundary rectifiability theorem, see \cite[Sec.~30]{Si}, we infer that $\rP(u)$ is i.m. rectifiable in $\R_{n-2}(B^n)$.
Furthermore, we already know that $\rP(u)$ has no inner boundary, see \eqref{bdPu}.
\par
Setting now
$$ T=G_u+S_u $$
where $S_u$ is the $n$-current given by \eqref{Su}, it suffices to show that $T$ is a Cartesian current. Since in fact $S_u$ is completely vertical, by uniqueness of the optimal lifting Cartesian current we readily obtain that $T=T_u$.
By the structure theorem, see \cite[Ch.~4]{GMSl1}, since we have just obtained that $S_u$ is i.m. rectifiable in $\R_n(B^n\tim\gR^2)$,
it suffices to show that $T$ satisfies the null-boundary condition
\beq\label{bdT} (\pa T)\pri B^n\tim\gR^2=0 \,. \eeq
In fact, we have:
$$ (\pa T)\pri B^n\tim\gR^2=(\pa G_u)\pri B^n\tim\gR^2+(-1)^{n-2}(\pa (\rP(u)\tim\qu{D^2}))\pri B^n\tim\gR^2 $$
where by the definition of boundary of a product of currents
$$
(\pa (\rP(u)\tim\qu{D^2}))\pri B^n\tim\gR^2=(\pa \rP(u))\pri B^n\tim\qu{D^2}+(-1)^{n-2}\rP(u)\tim\pa\qu{D^2}$$
so that \eqref{bdT} follows from \eqref{bdGucont}, \eqref{bdPu}, and property $\pa\qu{D^2}=\qu{\Sph^1}$.
\epf
\br As a consequence, in high dimension $n\geq 3$, by the previous result we infer that the distributional Jacobian $J(u)$ can be viewed as an $\gR^{d(n)}$-valued measure, with $d(n)=n(n-1)/2$, that is concentrated on the $(n-2)$-rectifiable set of points of positive multiplicity of the current $\rP(u)$, and actually
$$ |J(u)|(B^n)=\pi\cdot\gM(\rP(u))<\infty\,. $$
\er
\section{Examples}\label{Sec:example}
In this section, we give an easier example showing the strategy in our proof. We then show the existence of Sobolev maps in $W^{1,1}(B^n,\Sph^1)$ for which the relaxed energy is not finite.
\subsection{A model example}
Let $u\in W^{1,1}(B^n,\Sph^1)$ be defined as $$u(x)=\frac{\wid x}{|\wid x|}\,,\qquad x=(\wid x, \wih x)\in \gR^2\tim\gR^{n-2}.$$
If $n\geq 3$, the singular set of $u$ is the $(n-2)$-disk
$$ \DD^{n-2}:=\{(0,0,\wih x)\in\gR^n\,:\, |\wih x|\leq 1\}. $$
With the previous notation we get $\rP({u})=(-1)^{n-2}\qu{\DD^{n-2}}$, so that
by \eqref{bdGucont}
$$ (\pa G_{u})\pri B^n\tim\gR^2=(-1)^{n-1}\qu{\DD^{n-2}}\tim\qu{\Sph^1}. $$
This way we can equivalently write the lower bound for the relaxed energy as
$$ \int_{B^n}\sqrt{1+|\nabla u|^2}\,dx+|J(u)|(B^n)$$
where in low dimension $n=2$ we clearly have $J(u)=\Det \nabla u$.

\begin{figure}[t]
    \centering
    \includegraphics[scale=0.4]{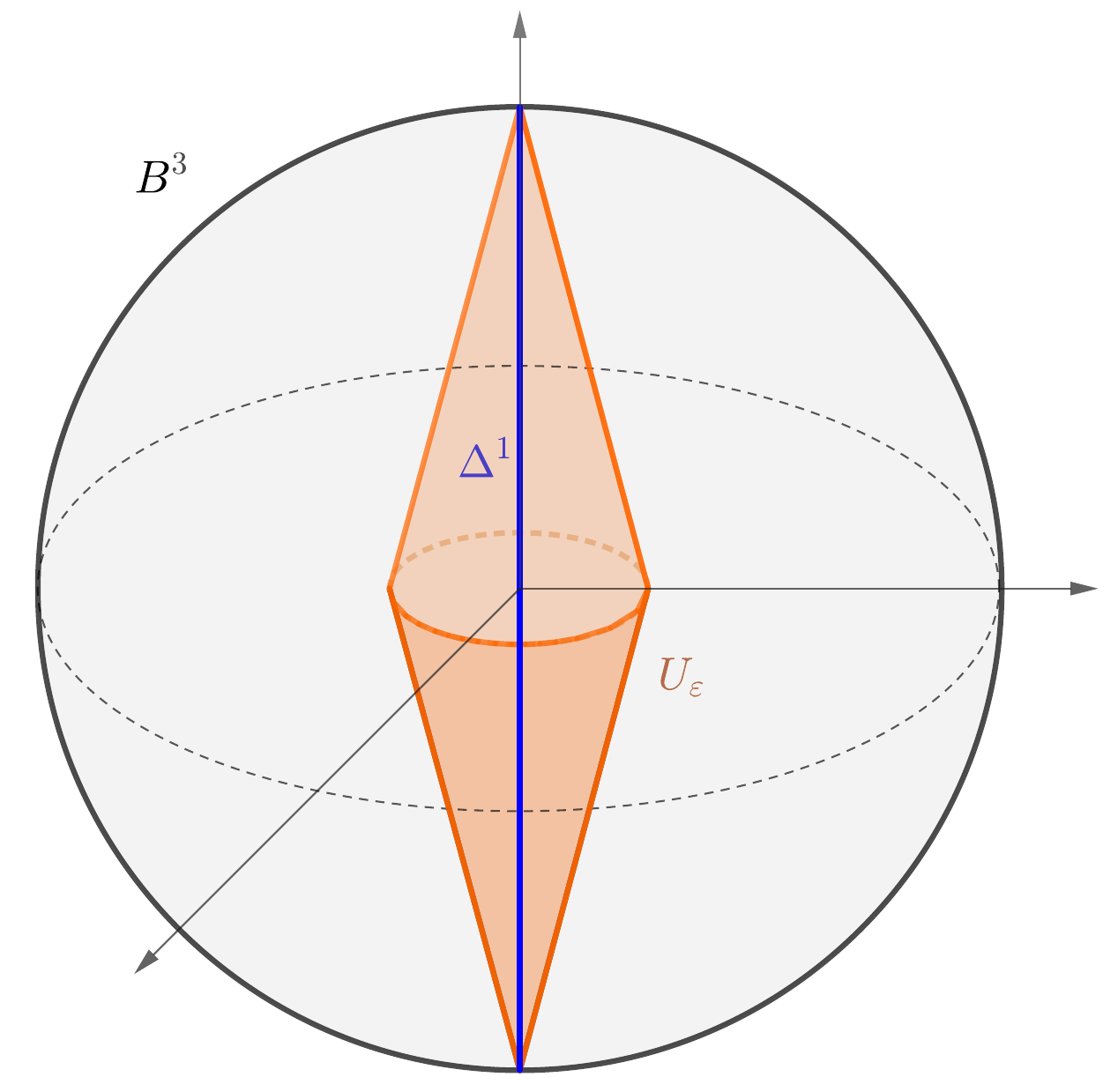}
      \caption{The cone shaped neighborhood $U_\e$ depicted in dimension $n=3$.}
    \label{double cone}
\end{figure}

For the upper bound estimate, we define a recovery sequence by constructing for each $\e>0$ small a suitable cone shaped neighborhood $U_\e$ of $\DD^{n-2}$ in the following way:
$$
U_\e:=\{x\in B^n: |\wid x|\leq\e(1-|\wih x|)\}
$$
and by defining $u_\e\in C^1(B^n,\gR^2)$ as
\begin{equation}
u_\e(x):=\left\{
\begin{aligned}
&u(x)\quad&&\mbox{ if }x\in B^n\setminus U_\e,\\
&\frac{|\wid x|}{\e(1-|\wih x|)}u(x)&&\mbox{ if }x\in U_\e.
\end{aligned}
\right.
\end{equation}
In the case $n=3$, $U_\e$ is a (double) cone of basis the disk $\wid{B}_\e:=\{x\in B^3: x=(\wid x,0), |\wid x|\leq\e\}$ and of vertices the North and South Poles of $B^3$ (see Fig. \ref{double cone}).\\
Let us check that $u_\e\to u$ in $W^{1,1}(B^n,\gR^2)$ as $\e\to0$. Clearly, $\int_{U_\e}|u_\e|\,dx\to0$, since $|u_\e|\leq|u|=1$.
Therefore, it is enough to prove that
\beq\label{conv}
\lim_{\e\to 0}\int_{U_\e}|\nabla u_\e(x)|dx=0\,.
\eeq

In cylindrical coordinates
$$
\bar u_\e(\rho,\theta,\wih x):=u_\e(\rho\cos\theta,\rho\sin\theta,\wih x),\,\, \rho\in[0,1],\,\,\theta\in[0,2\pi),\,\, \wih x\in\gR^{n-2}\,, $$
where $|\wih x|\leq1$,
we have
$$
\bar u_\e(\rho,\theta,\wih x)=\frac{\rho}{\e(1-|\wih x|)}(\cos\theta,\sin\theta)\quad{\textrm{if }}\rho\leq \e(1-|\wih x|)\,.
$$
Compute the partial derivatives of $\bar u_\e$:
\begin{align*}
&\partial_\rho \bar u_\e(\rho,\theta,\wih x)=\frac{1}{\e (1-|\wih x|)}(\cos\theta,\sin\theta),\\
&\partial_\theta \bar u_\e(\rho,\theta,\wih x)=\frac{\rho}{\e (1-|\wih x|)}(-\sin\theta,\cos\theta),\\
& \partial_{\wih x}\bar u_\e(\rho,\theta,\wih x)=\frac{\rho}{\e (1-|\wih x|)^2}(\cos\theta,\sin\theta)\otimes\frac{\wih x}{|\wih x|}.
\end{align*}
Moreover, by identifying the set $\{\wih x\in\gR^{n-2}:|\wih x|\leq 1\}$ with $\DD^{n-2}$, we have 
$$
\begin{array}{l}
\ds \int_{U_\e}|\nabla u_\e(x)|dx = \\
\qquad \ds =\int_{\DD^{n-2}}\int_0^{2\pi}\int_{0}^{\e(1-|\wih x|)}\rho\,\sqrt{|\partial_\rho \bar u_\e|^2+\frac{|\partial_\theta \bar u_\e|^2}{\rho^2}+|\partial_{\wih x}\bar u_\e|^2}\,d\rho d\theta d\wih x\\
\qquad \ds =\int_{\DD^{n-2}}\int_0^{2\pi}\int_{0}^{\e(1-|\wih x|)}\rho\,\sqrt{\frac{2}{\e^2(1-|\wih x|)^2}+\frac{\rho^2}{\e^2(1-|\wih x|)^4}}\,d\rho d\theta d\wih x\\
\qquad\ds \leq\int_{\DD^{n-2}}\int_0^{2\pi}\int_{0}^{\e(1-|\wih x|)}\left[\frac{2\rho}{\e(1-|\wih x|)}+\frac{\rho^2}{\e(1-|\wih x|)^2}\right]\,d\rho d\theta d\wih x\\
\qquad\ds \leq\int_{\DD^{n-2}}\int_0^{2\pi}\int_{0}^{\e(1-|\wih x|)}\left[2+\e\right]\,d\rho d\theta d\wih x
\to 0\quad\mbox{ as }\e\to 0^+,
\end{array}
$$
where we used that $\rho=|\wid x|\leq\e(1-|\wih x|)$ in $U_\e$.
Therefore, \eqref{conv} holds, and by dominated convergence
\begin{align}\label{eq:dom conv}
\lim_{\e\to0}\int_{B^n}\sqrt{1+|\nabla u_\e |^2}\,dx=\int_{B^n}\sqrt{1+|\nabla u|^2}\,dx\,.
\end{align}
It remains to check that
$$ \limsup_{\e\to0}\int_{B^n}|M_2(\nabla u_\e)|dx\leq\pi\,\mathcal{H}^{n-2}({\DD^{n-2}})=\pi\,\gM(\rP(u))\,. $$
We have $|M_2(\nabla u_\e)|=|M_2(\nabla \bar u_\e)|$, where we compute the components of $M_2(\nabla \bar u_\e)$ w.r.t. the basis in cylindrical coordinates:
\begin{align*}
&M_2(\nabla \bar u_\e)_{12}=\frac 1\rho\,\partial_\rho \bar u_\e\wedge\partial_\theta \bar u_\e=\frac{1}{\e^2(1-|\wih x|)^2}\,,\\
&M_2(\nabla \bar u_\e)_{1j}=\partial_\rho \bar u_\e\wedge\partial_{x_j} \bar u_\e=0\quad&\forall j=3,\ldots,n\,, \\
&M_2(\nabla \bar u_\e)_{2j}=\frac 1\rho\partial_\theta \bar u_\e\wedge\partial_{x_j} \bar u_\e=-\frac{\rho}{\e^2(1-|\wih x|)^3}\frac{x_j}{|\wih x|}\quad&\forall\, j=3,\ldots,n\,, \\
&M_2(\nabla \bar u_\e)_{ij}=\partial_{x_i} \bar u_\e\wedge\partial_{x_j} \bar u_\e=0\quad&\forall\, i,j=3,\ldots,n,\,\, i\neq j\,.
\end{align*}
Therefore,
$$\begin{array}{l}
\ds \int_{B^n}|M_2(\nabla u_\e)|dx = \\
\qquad\ds =\int_{U_\e}|M_2(\nabla u_\e)|dx=\int_{\DD^{n-2}}\int_0^{2\pi}\int_{0}^{\e(1-|\wih x|)}\rho\,|M_2(\nabla\bar u_\e)|\,d\rho d\theta d\wih x\\
\qquad\ds \leq\int_{\DD^{n-2}}\int_0^{2\pi}\int_{0}^{\e(1-|\wih x|)}\left[\frac{\rho}{\e^2(1-|\wih x|)^2}+\frac{\rho^2}{\e^2(1-|\wih x|)^3}\right]\,d\rho d\theta d\wih x\\
\qquad\ds =\int_{\DD^{n-2}}\int_0^{2\pi}\frac{1}{2}d\theta d\wih x+\int_{\DD^{n-2}}\int_0^{2\pi}\int_{0}^{\e(1-|\wih x|)}\frac{\rho^2}{\e^2(1-|\wih x|)^3}\,d\rho d\theta d\wih x\\
\qquad\ds =\pi\,\mathcal{H}^{n-2}(\DD^{n-2})+O(\e)
\to\pi\,\mathcal{H}^{n-2}(\DD^{n-2})\quad\mbox{ as }\e\to0^+.
\end{array}
$$
Using \eqref{eq:dom conv}, we conclude
\begin{align*}
\limsup_{\e\to0}\gM(G_{u_\e})&\leq\lim_{\e\to0}\int_{B^n}\sqrt{1+|\nabla u_\e |^2}dx+\limsup_{\e\to0}\int_{B^n}|M_2(\nabla u_\e)|dx\\
&\leq\int_{B^n}\sqrt{1+|\nabla u|^2}dx+\pi\mathcal{H}^{n-2}(\DD^{n-2})\,.
\end{align*}
\br
In the previous example, the choice of the cone shaped neighborhood is not crucial for the computation of the upper bound estimate. We could have followed essentially the same argument also by taking $U_\e$ of cylindrical shape, i.e. by defining $U_\e:=(\wid{B}_\e\times\DD^{n-2})\cap B^n$. The advantage of the cone shaped construction is that the width of $U_\e$ shrinks at the boundary of $\DD^{n-2}$, which will be useful in the case the singular set of $u$ is polyhedral.
\er
\subsection{Sobolev maps with unbounded relaxed energy}
We show the existence of Sobolev maps $u\in W^{1,1}(B^n,\Sph^1)$ which do not have finite relaxed energy.

In low dimension $n=2$, it suffices to find a sequence $\{B_j\}$
of pairwise disjoint balls contained in $B^2$ such that the
restriction $u_{\vert B_j}$ behaves like a vortex map around the center of $B_j$.
Therefore, by the superadditivity of the
set function corresponding to the localization of the relaxed energy, we obtain a contribution equal to $\pi$ around each singular point.
In particular, $|\Det \nabla u|(B^2)=\infty$.

The counterexample in high dimension $n\geq 3$ is trivially obtained by setting $\ol u(x)=u(x_1,x_2)$, for $x\in B^n$.
In that case, we clearly have $|J(u)|(B^n)=\infty$.
\smallskip
\par
\begin{figure}[t]
    \centering
    \includegraphics[scale=0.4]{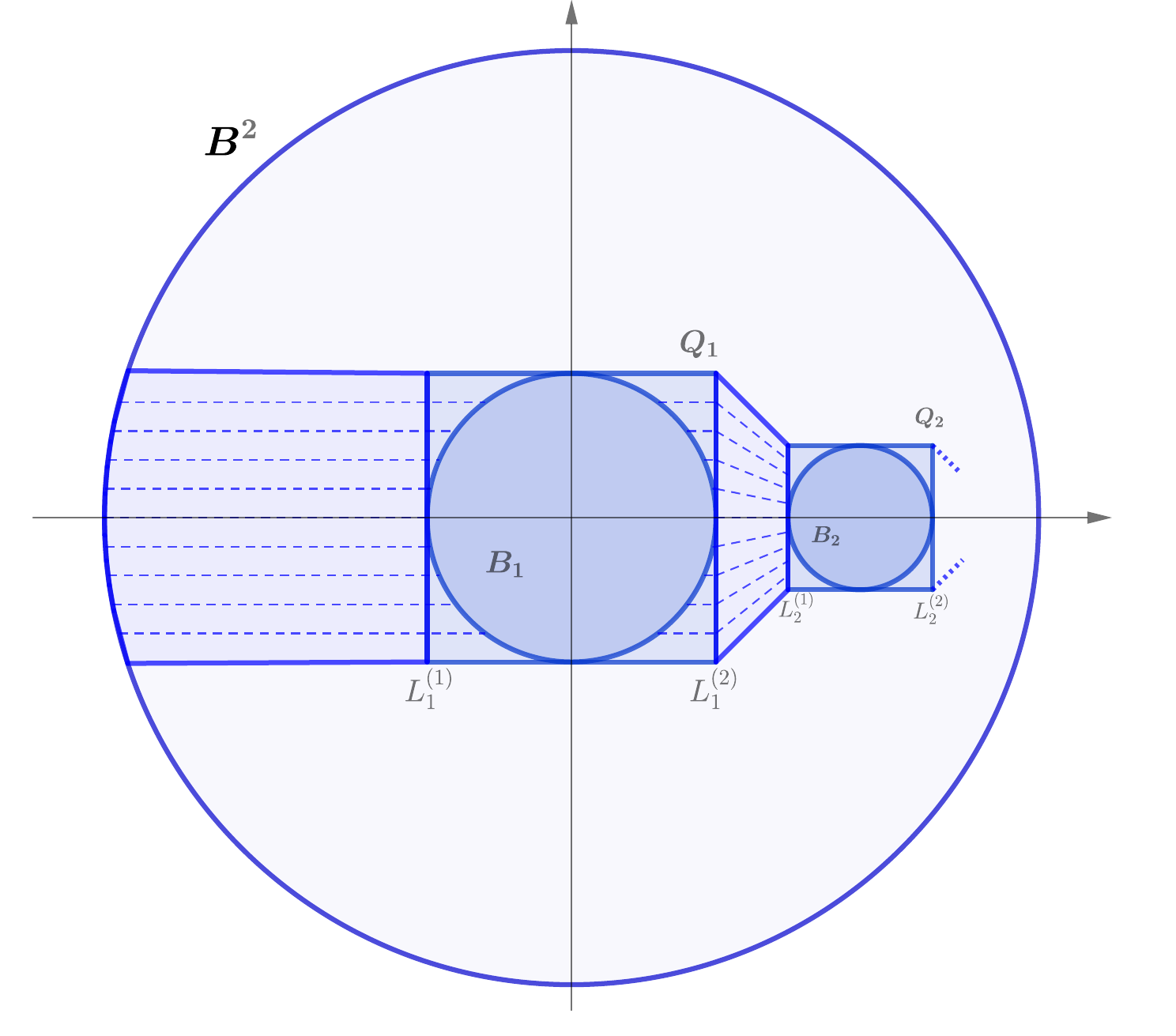}
      \caption{The construction in the source disk $B^2$. On each disk $B_j$ the vortex map is replicated with alternating orientation. }
    \label{domain}
\end{figure}

 Following an example by \cite{Mu06}, we set
$B_j:=B^2(c_j,2^{-(j+1)})$, where
$$
c_j=(1-2^{1-j},0)\,,\qquad j=1,2,\ldots
$$
Moreover we define $u_{\vert B_j}:=u^{(j)}:B_j\to\gR^2$ by
$$
u^{(j)}(x):=\left\{
\renewcommand{\arraystretch}{1.8} \ba{ll}
\ds\frac{x-c_j}{\vert x-c_j\vert} & \If\quad
j=1,3,5,\ldots   \\
\q\biggl(\ds\frac{x-c_j}{\vert x-c_j\vert}\biggr) & \If\quad
j=2,4,6,\ldots  \ea\right.
$$
where $\q:\Sph^1\to\Sph^1$ is
defined (in terms of the angle function $\theta$ on $\Sph^1$) by
$$
\q(\theta):=-\theta+\p\,.
$$
If $Q_j:=c_j+[-2^{-(j+1)},2^{-(j+1)}]^{\,2}$ denotes the
square circumscribing $B_j$, we extend $u_{\vert B_j}$ to
$Q_j$ as the continuous map which is constant in the
$x_1$-variable (note that $Q_j\subset B^2$ for every $j\geq 1$, see Fig. \ref{domain}).
Then $u\equiv (1,0)$ and $u\equiv (-1,0)$ over all the upper and respectively lower sides of the
boundary of the $Q_j$'s which are parallel to the $x_1$-axis,
whereas on the sides parallel to the $x_2$-axis,
$$
L^k_j:=c_j+\{((-1)^k\,2^{-(j+1)},x_2)\mid -2^{-(j+1)}\leq x_2\leq
2^{-(j+1)}\}\,, \qquad k=1,2\,,
$$
both $u_{\vert L^2_j}$ and $u_{\vert L^1_{j+1}}$ parameterize
the same half of the circle $\Sph^1$ with the same orientation. We can thus define
$u$ over the convex hull of $L^2_j$ and $L^1_{j+1}$, the
right-hand side of $\pa Q_j$ and the left-hand side of $\pa
Q_{j+1}$, as the continuous map which is constant along the
straight lines connecting the corresponding points in $L^2_j$
and $L^1_{j+1}$ (points on which $u$ takes the same
value). We finally define $u$ in the strip connecting
$L^1_1$ to the boundary of $B^2$ as the continuous map
constant in the $x_1$-variable, and set $u\equiv (1,0)$ or $u\equiv (-1,0)$ in the
two remaining components of $B^2$. Then, it is not difficult to show that $u\in
W^{1,1}(B^2,\gR^2)$. On the other hand, by the result from \cite{BSS} we know that for each $j$, the relaxed energy of $u^{(j)}$ on $B_j$ is greater than $\p$.
Therefore, by the superadditivity of the (localized) relaxed functional it turns out that the map $u$ does not have a finite relaxed energy.
\section{The explicit formula}\label{Sec:main}
The Main Result of this paper is the following
\bt\label{Tmain} Let $n\geq 2$ and $u\in W^{1,1}(B^2,\Sph^1)$. Then, $\ol{\mathcal{A}}_{\BV}(u)<\i$ if and only if the $(n-2)$-current $\rP(u)$ is i.m. rectifiable and with finite mass, $\gM(\rP(u))<\i$.
In that case, moreover, one has:
$$ \ol{\mathcal{A}}_{\BV}(u) =\int_{B^n}\sqrt{1+|\nabla u|^2}\,dx+\pi\,\gM(\rP(u))\,. $$
\et

In dimension $n=2$, recalling that $\pi\,\gM(\rP(u))=|\Det\nabla u|(B^n)$, the latter result was proved in \cite{BSS}.
In high dimension $n\geq 3$, the {\em energy gap}, $\pi\,\gM(\rP(u))$, agrees with the total variation of the distributional Jacobian $J(u)$.
\subsection{Energy lower bound}
By the previous results, we readily obtain the energy lower bound:
\bp\label{Plbd} If $u\in W^{1,1}(B^n,\Sph^1)$ has finite relaxed energy \eqref{Area-rel}, then
$$ \ol{\mathcal{A}}_{\BV}(u)\geq\int_{B^n}\sqrt{1+|\nabla u|^2}\,dx+\pi\,\gM(\rP(u))\,. $$
\ep
\bpf Choose any smooth sequence $\{u_h\}\sb C^\infty(B^n,\gR^2)$ such that $u_h\to u$ in $L^1(B^n,\gR^2)$ and $\int_{B^n}|\nabla u_h|\,dx\to\int_{B^n}|\nabla u|\,dx$ as $h\to\infty$, and such that
$$\sup_h\int_{B^2}|M_2(\nabla u_h)|\,dx<\infty\,. $$
Possibly taking a not relabeled subsequence, we may and do assume that
$$ \liminf_{h\to\i}\mathcal{A}(u_h)=\lim_{h\to\i}\mathcal{A}(u_h)<\infty\,. $$
Since $(\pa G_{u_h})\pri B^n\tim\gR^2=0$ and the mass of $G_{u_h}$ is given by \eqref{MGu}, with $u=u_h$,
by applying Federer-Fleming's closure theorem, see \cite[Sec.~32]{Si}, and on account of the strict convergence $u_h\BVc u$, it turns out that possibly passing to a subsequence, $G_{u_h}\wc T_u$ weakly in $\D_n(B^n\tim\gR^2)$ to the unique optimal lifting Cartesian current $T_u$, so that by lower semicontinuity of the mass
$$ \gM(T_u)\leq\liminf_{h\to\infty}\gM(G_{u_h})=\lim_{h\to\i}\mathcal{A}(u_h)\,. $$
Since we already know that
$$\gM(T_u)=\gM(G_u)+\gM(S_u)=\int_{B^n}\sqrt{1+|\nabla u|^2}\,dx+ \pi\,\gM(\rP(u)) $$
the energy lower bound readily follows.
\epf
\subsection{The approximation theorem}
The energy upper bound, which yields to the validity of Theorem~\ref{Tmain}, is an immediate consequence of the following approximation result:
\bt\label{Tappr} Let $n\geq 2$ and $u\in W^{1,1}(B^n,\Sph^1)$ be a Sobolev map with finite relaxed energy \eqref{Area-rel}. Then, there exists a smooth sequence
$\{u_h\}\sb C^\infty(B^n,\gR^2)$ such that $G_{u_h}\wc T_u$ weakly in $\D_n(B^n\times\gR^2)$ and
$\gM(G_{u_h})\to \gM(T_u)$ as $h\to\infty$.
\et
\par In fact, the weak convergence with the mass implies the strict $\BV$-convergence\footnote {The strata of $G_{u_h}$ must converge to the corresponding ones of $T_u$ in measure and in total variation.} and the energy limit
$$ \lim_{h\to\i}\mathcal{A}(u_h)=\int_{B^n}\sqrt{1+|\nabla u|^2}\,dx+\p\,\gM(\rP(u))\,. $$
Therefore, if $\ol{\mathcal{A}}_{\BV}(u)<\i$, by the explicit formula we obtain that $\gM(\rP(u))<\i$. On the other hand, when $\gM(\rP(u))<\i$, the approximation theorem~\ref{Tappr} continues to hold, yielding to the optimal upper bound and hence to condition $\ol{\mathcal{A}}_{\BV}(u)<\i$. Therefore,
Theorem~\ref{Tmain} holds true.
\smallskip\par In low dimension, the approximating sequence is readily obtained:
\smallskip\par\noindent
\bpf[Proof of Theorem \ref{Tappr}, case $n=2$.]
By Bethuel's results in \cite{Bet}, we can find a sequence $\{u_h\}\sb W^{1,1}(B^2,\Sph^1)$ strongly converging to $u$ in $W^{1,1}(B^2,\gR^2)$ and such that each $u_h$ is smooth outside a finite set of points. Furthermore, we have
\beq\label{MGPuh}\lim_{h\to\i}\gM(\rP(u_h))=\gM(\rP(u))\,. \eeq
\par In fact, for any square $F$ of the grid in Bethuel's proof, the restriction of $u$ to the boundary of $F$ is a continuous function with Brouwer degree $d_F\in\gZ$ satisfying $|d_F|\leq \gM(\rP(u)\pri \Int(F))$. As a consequence, it turns out that $\gM(\rP(u_h))\leq\gM(\rP(u))$ for each $h$.
Therefore, by lower semicontinuity we obtain \eqref{MGPuh}.
\par We now show that for each $h$ we can find a smooth sequence $\{u^{(h)}_k\}$ in $C^\i(B^2,\gR^2)$ strongly converging to $u$ in $W^{1,1}(B^2,\gR^2)$ and such that
$$ \lim_{k\to\i}\mathcal{A}(u^{(h)}_k)=\int_{B^2}\sqrt{1+|\nabla u_h|^2}\,dx+\pi\,\gM(\rP(u_h))\,. $$
\par Since we make use of a local argument, without loss of generality we may and do assume that $v=u_h$ is smooth outside the origin and $\rP(v)=d\,\d_{0_{\gR^2}}$ for some $d\in\gZ$.
\par
For every $\e>0$ small, the restriction $v_{\vert\pa B^2_\e}$ is a smooth map of degree $d$. Therefore,
we can find a smooth homotopy map $H:[0,1]\tim[0,2\pi]\to\Sph^1$ such that
$H(0,\t)=(\cos(d\theta),\,\sin(d\theta))$ and $H(1,\t)=v(\e\cos\t,\e\sin\t)$, where we have introduced the standard polar coordinates $x=\r(\cos\t,\sin\t)$.
Define now $v_\e:B^2_\e\to\gR^2$ as
$$v_\e(\r\cos\t,\r\sin\t):=\left\{ \ba{ll} H\Bigl( 2\r/\e-1,\t\Bigr) & {\text{if }}
\e/2\leq \r\leq \e \\
(2\r/\e)\,\bigl(\cos(d\t),\sin(d\t)\bigr) & {\text{if }}
 \r\leq \e/2 \,.
\ea \right. $$
It is readily checked that
$$\mathcal{A}(v_\e)\leq \int_{B^2}\sqrt{1+|\nabla v|^2}\,dx+\pi\,|d|+O(\e) $$
where $O(\e)\to 0$ as $\e\to 0$. Since moreover the graph currents $G_{v_\e}$ weakly converge to the Cartesian current $T:=G_v+d\,\d_0\tim\qu{D^2}$ along a sequence $\e_h\to 0$, by lower semicontinuity we obtain
$$ \lim_{h\to\i}\mathcal{A}(v_{\e_h})=\lim_{h\to\i}\gM(G_{v_{\e_h}})=\gM(T)=\int_{B^2}\sqrt{1+|\nabla v|^2}\,dx+\pi\,|d| $$
where $|d|=\gM(\rP(u))$.
Further details are omitted. \epf
%
\par In the sequel we therefore assume $n\geq 3$. Theorem~\ref{Tappr} is obtained by applying the following technical results,
the proof of which is collected in the next section.
\subsection{Reduction to maps with a nice singular set}
Firstly, we find an approximating sequence which is smooth outside a singular set given by (the support of) a polyhedral chain, in such a way that we have mass convergence of the current of the singularities.
\par Let $Q^n=]-1,1[^n$ denote the open $n$-cube in $\gR^n$ of side two, and let $u\in W^{1,1}(Q^n,\Sph^1)$.
Denote by $R^\i(Q^n,\Sph^1)$ the subclass of maps $u$ in $W^{1,1}(Q^n,\Sph^1)$ which are smooth outside a nice singular set $\sing u$ of codimension two.
This means that $\sing u$ is given by the support of some polyhedral $(n-2)$-chain $\rP$ in $Q^n$, and actually $\rP(u)=\rP$.
More precisely, we have
\beq\label{Pol} \rP=\sum_{i=1}^m d_i\qu{\Delta_i}\,,\quad \gM(\rP)=\sum_{i=1}^m |d_i|\,\Ha^{n-2}({\Delta_i})<\infty \eeq
for some $m\in\Nat^+$, where $d_i\in\gZ$ and $\Delta_i$ is an oriented $(n-2)$-simplex contained in the closure of $Q^n$, for each $i$. Notice that $d_i$ coincides with the degree of $u$ around $\DD_i$ up to a sign, precisely
$\deg(u,\DD_i)=(-1)^{n-2}d_i$.
The support $\spt\rP$ of $\rP$ is the union of the closures of the simplices $\Delta_i$. Moreover, after a subdivision we may and do assume that
$\Int(\DD_i)\cap\Int(\DD_j)=\emp$ for $1\leq i<j\leq m$, so that the simplices $\DD_i$ and $\DD_j$ possibly intersect at points in the common $(n-3)$-skeleton.
\par By Bethuel's theorem, the class $R^\i(Q^n,\Sph^1)$ is dense in $W^{1,1}(Q^n,\Sph^1)$ strongly in $W^{1,1}(Q^n,\gR^2)$.
To our purposes, we shall see that for maps with finite relaxed energy, something more can be said.
\bt\label{TBet} Assume that $u\in W^{1,1}(Q^n,\Sph^1)$ has finite relaxed energy.
Then, we can find a sequence $\{u_k\}\sb R^\i(Q^n,\Sph^1)$ strongly converging to $u$ in $W^{1,1}(Q^n,\gR^2)$ and such that
$$ \lim_{k\to\i}\gM(\rP(u_k))=\gM(\rP(u))\,. $$
\et
\subsection{Energy approximation at the singular set}
Assume now that $n\geq 3$ and $u\in R^\i(Q^n,\Sph^1)$ satisfies $\gM(\rP(u))<\i$, with $\rP(u)=\rP$ as in \eqref{Pol}.
Without loss of generality, we assume that for e.g. $i=1$
we have
$\DD_1=\{0_{\gR^2}\}\tim\wih \DD$. 
Let
$$ \DD_\e:=\{(\wid x,\wih x)\in\gR^2\tim\wih\DD \,:\, \vert\wid x\vert\leq \e y(\wih x) \}\,,\quad \e>0 $$
where we have denoted
\beq\label{dist}
y(\wih x):=\dist(\wih x,\pa\wih\DD)\,,
\eeq
so that for $\e>0$ small, the cone $\DD_\e$ intersects the other simplices $\DD_i$
only at points in $\pa\DD_i$, for $i=2,\ldots,m$.
Since moreover $u\in W^{1,1}(Q^n,\Sph^1)$ is smooth outside the support of $\rP(u)$, for a.e. $\e>0$ the restriction of $u$ to the boundary of $\DD_\e$ is in $W^{1,1}$.
Furthermore, recalling the definition of the current $\rP(u)$, it turns out that for all $\wih x\in \wih\DD$ the degree of $u(\cdot,\wih x)$ around
$0_{\gR^2}$ is constantly equal to $d=d_i\in\gZ$.
\par The following result allows to remove the dipole $\Delta:=\DD_1$, by paying an amount of energy essentially equal to $\pi\,|d|\,\Ha^{n-2}(\Delta)$. Of course, the argument will be applied to each simplex $\DD_i$ in the proof of Theorem~\ref{Tappr}.
\bt\label{Tdipolo} For a.e. $\e>0$ small, there exists a smooth map $v_\e:\DD_\e \to\gR^2$ such that $v_{\e\vert\pa\DD_\e}=u_{\vert\pa\DD_\e}$ in the sense of
the traces, and
\beq\label{est-dipolo} \gM(G_{v_\e}\pri\DD_\e\tim\gR^2)\leq \pi\,|d|\,\Ha^{n-2}(\DD)+O(\e)\,.
\eeq
\et
\subsection{Removal of point singularities}
In dimension $n=3$, we also need the following argument that allows to remove point singularities by paying a small amount of energy.
\bt\label{Tsing} Let $n\geq 3$ and $u\in W^{1,1}(Q^n,\gR^2)$ be smooth outside a discrete set.
Then there exists a sequence $\{u_k\}\sb C^\i(B^n,\gR^2)$ such that $u_k\BVc u$ strictly, and
$$ \lim_{k\to\infty}\int_{Q^n}\sqrt{1+|\nabla u_k|^2+|M_2(\nabla u_k)|^2}\,dx=\int_{Q^n}\sqrt{1+|\nabla u|^2+|M_2(\nabla u)|^2}\,dx\,. $$
\et
\subsection{Removal of high codimension singularities}
In high dimension $n\geq 4$, instead, we first have to remove singularities of codimension greater than two.
More precisely, let $k=3,\ldots n-1$, integer, and let $\Delta$ denote an $(n-k)$-dimensional simplex contained in the closure of $Q^n$.
Without loss of generality, assume $ \DD=\{0_{\gR^k}\}\tim\wih \DD$.
For $\e>0$ small, define again
$$ \DD_\e:=\{(\wid x,\wih x)\in\gR^k\tim\wih\DD \,:\, |\wid x|\leq \e y(\wih x) \} $$
where $y(\wih x)$ is the distance function in \eqref{dist}.
\bt\label{Tsimplex} Let $u\in W^{1,1}(B^n,\gR^2)$ be smooth in $\Int(\DD_{\e_0})\sm \DD$ for some $\e_0>0$ small.
Then, for a.e. $\e>0$ small, there exists a smooth map $v_\e:\Int(\DD_\e) \to\gR^2$ such that $v_{\e\vert\pa\DD_\e}=u_{\vert\pa\DD_\e}$ in the sense of
the traces, and
\beq\label{est-simplex} \gM(G_{v_\e}\pri\DD_\e\tim\gR^2)\leq O(\e)\,.
\eeq
\et
\subsection{Proof of the approximation theorem}
We are now in position to give the:
\smallskip\par\noindent
\bpf[Proof of Theorem \ref{Tappr}, case $n\geq 3$.] Since the weak convergence of currents with supports contained in the closure of $B^n\tim D^2$ is metrizable, compare \cite[Sec.~31]{Si}, we can apply a diagonal argument. Moreover, since $B^n$ is bilipschitz homeomorphic to $Q^n$, we may and do assume $u:Q^n\to\Sph^1$.

\noindent
{\em Step 1.} By Theorem~\ref{TBet}, we reduce to the case in which $u\in R^\i(Q^n,\Sph^1)$ and $u$ is smooth outside the support of $\rP(u)$, a polyhedral $(n-2)$-chain in $Q^n$.

\noindent
{\em Step 2.}
By applying iteratively Theorem~\ref{Tdipolo}, for each $\e>0$ we find a Sobolev map $u_\e\in W^{1,1}(Q^n,\gR^2)$ that is smooth outside an $(n-3)$-dimensional polyhedral chain $\SS_\e$ and such that
$$ \int_{Q^n}\sqrt{1+|\nabla u_\e|^2+|M_2(\nabla u_\e)|^2}\,dx\leq \int_{Q^n}\sqrt{1+|\nabla u|^2}\,dx+\pi\,\gM(\rP(u))+\e $$
with $u_\e\BVc u$ strictly and $\gM(\rP(u_\e))\to\gM(\rP(u))$, as $\e\to 0$.

\noindent
{\em Step 3.} If $n=3$, the finite set $\SS_\e$ of point singularities of $u_\e$ is removed by means of Theorem~\ref{Tsing}.
In high dimension $n\geq 4$, we first apply Theorem~\ref{Tsimplex}, with $k=3$, to each $(n-3)$-simplex of $\SS_\e$, and reduce to the case of a map
that is smooth outside an $(n-4)$-dimensional polyhedral chain, given by the union of the faces $F$ of the $(n-3)$-simplices of $\SS_\e$ that lie inside $Q^n$. If $n\geq 5$, we then iteratively repeat the same argument, by applying Theorem~\ref{Tsimplex} for $k=4,\ldots,n-1$. Finally, we apply Theorem~\ref{Tsing} in order to remove the finite set of point singularities.

\noindent
{\em Step 4.} By a diagonal argument, we find a good approximating sequence given by Lipschitz-continuous functions, where the weak convergence as currents readily follows. By a standard convolution argument, the proof is complete.
\epf
\section{Proofs}\label{Sec:proofs}
In this section we collect the proofs of the technical results leading to Theorem~\ref{Tappr}.
Recall that we assume $n\geq 3$.
\subsection{Reduction to maps with a nice singular set}
\bpf[Proof of Theorem \ref{TBet}.]
We first consider the case when $u$ is smooth in a neighborhood of the boundary of $Q^n$.
Therefore, $\rP(u)$ can be viewed as a current in $\R_{n-2}(\gR^n)$ satisfying $\pa\rP(u)=0$, see Theorem~\ref{TPu}, and with support a closed set contained in the open cube $Q^n$. In the sequel, $c(n)$ will denote a real positive constant only depending on the dimension $n$, possibly varying from line to line.
\par
Let $\s>0$ small. By Federer's strong polyhedral approximation theorem, see \cite[4.2.2]{FH}, see also \cite[Sec.~2.2.6]{GMSl1}, there exists a diffeomorphism $\vf_\s$ of $Q^n$ onto itself and an
$(n-2)$-dimensional polyhedral chain $\rP_\s$ with support in $Q^n$, such that $\vf_{\s\#}\rP(u)-\rP_\s=\pa R_\s$ for some current
$R_\s\in\R_{n-1}(Q^n)$ with $\gM(R_\s)+\gM(\pa R_\s)<\s$.
Moreover, $\Lip\vf_\s\leq 1+\s$, $\Lip{\vf_\s}^{-1}<1+\s$, and $\vf_\s(x)=x$ if the distance of $x$ to the support of $\rP(u)$ is greater than $\s$.
%
\par
\begin{figure}[t]
    \centering
    \includegraphics[scale=0.35]{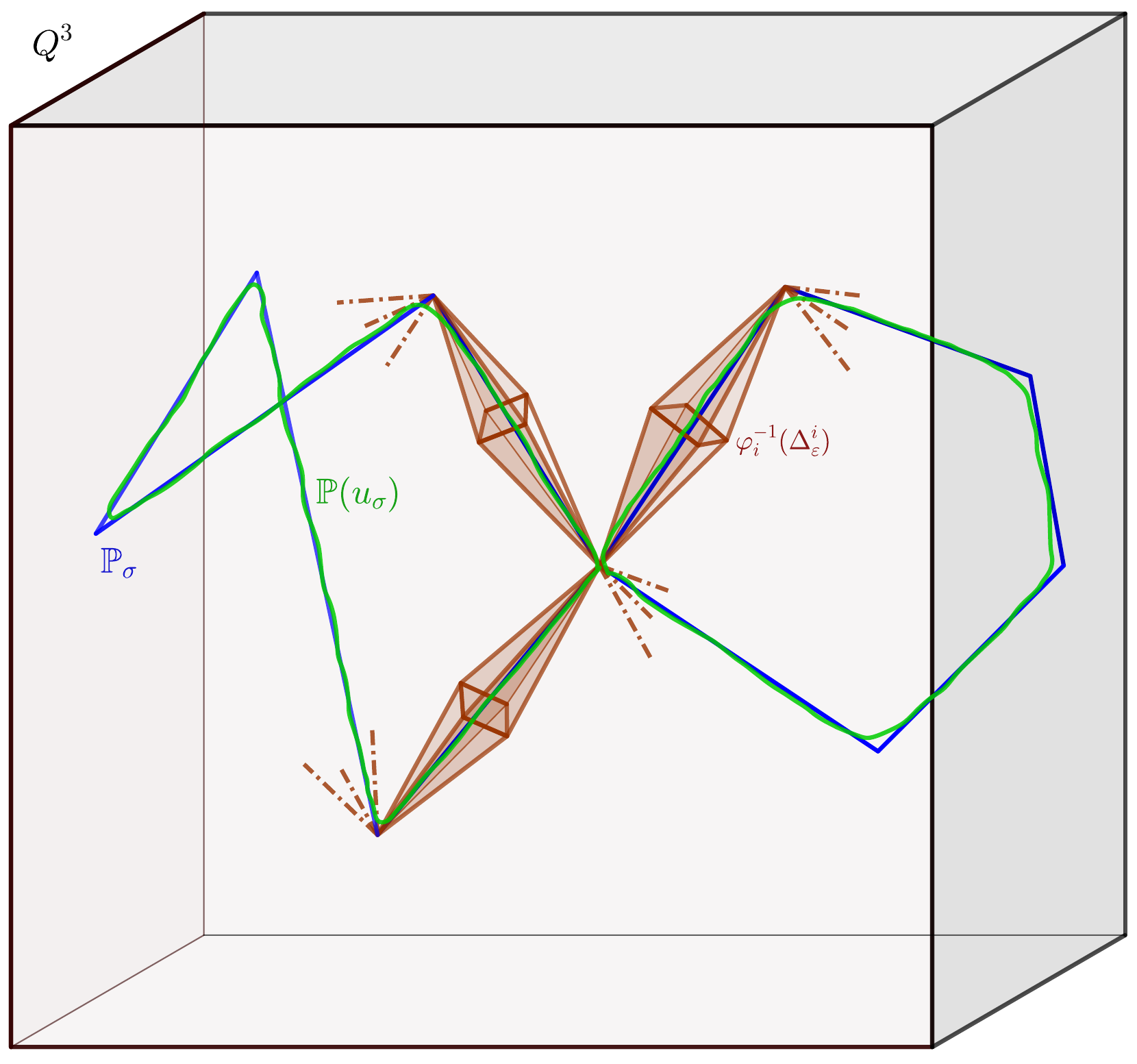}
      \caption{The polyhedral chain $\rP_\s$ (in blue), the current $\rP(u_\s)$ (in green) and the pyramidal neighborhoods ${\vf_i}^{-1}(\DD^i_\e)$, depicted in dimension $n=3$. Notice that both $\rP_\s$ and $\rP(u_\s)$ are boundaryless and have support contained in the open cube $Q^3$, since we are assuming that $u$ is smooth near the boundary of $Q^3$.}
    \label{pyramids}
\end{figure}

 Letting $u_\s:=u\circ{\vf_\s}^{-1}$, then $u_\s\in W^{1,1}(Q^n,\Sph^1)$, $u_\s\to u$ strongly in $W^{1,1}(Q^n,\gR^2)$ as $\s\to 0$, and
$ \rP(u_\s)=\vf_{\s\#}\rP(u)$ (see Fig. \ref{pyramids}), so that
\beq\label{Federer}
\rP(u_\s)-\rP_\s=\pa R_\s\,,\quad \gM(R_\s)+\gM(\pa R_\s)<\s\,.
\eeq
As a consequence, the open set
$$U_\s=Q^n\sm \spt \rP_\s $$
has full measure, and
\beq\label{sigma} \gM(\rP(u_\s)\pri U_\s)=\gM((\pa R_\s)\pri U_\s)\leq \gM(\pa R_\s)<\s\,. \eeq
\par For any $\s>0$ small, we now write $u=u_\s$, $\rP=\rP_\s$, and $U=U_\s$, for the sake of simplicity, and we write $\rP$ as in \eqref{Pol}.
After a rigid motion $\vf_i$ in $\gR^n$ we have
$$\vf_i(\DD_i)=\{0_{\gR^2}\}\tim\wih \DD_i \quad\fa\,i=1,\ldots,m\,. $$
If $\wid x=(x_1,x_2)\in\gR^2$, we let $\Vert\wid x\Vert:=|x_1|+|x_2|$, and for $\e>0$ small
$$ \DD^i_\e:=\{(\wid x,\wih x)\in\gR^2\tim\wih\DD_i \,:\, \Vert\wid x\Vert\leq \e y_i(\wih x) \} $$
where we have denoted $y_i(\wih x):=\dist(\wih x,\pa\wih\DD_i)$. Therefore, there exists $\e_0>0$ such that if $\e\in]0,\e_0[$, for any $1\leq i<j\leq m$ the cones $\DD^i_\e$ and $\DD^j_\e$ are
interiorly disjoint, and only intersect at points in the $(n-2)$-dimensional set $\DD^i\cap\DD^j$.
\par Denote by $\SS^i_\e(\ell)$ the $\ell$-dimensional skeleton of ${\vf_i}^{-1}(\DD^i_\e)$ (see Fig. \ref{pyramids}). By a slicing argument, it turns out that for a.e. $\e\in]0,\e_0[$ the restriction $u_{\vert F}$ of $u$ to any $\ell$-face $F$ of $\SS^i_\e(\ell)$ is a Sobolev map in $W^{1,1}(F,\Sph^1)$, for each $\ell=1,\ldots,n-1$ and $i=1,\ldots,m$.
In the sequel, we shall tacitly assume that $\e\in]0,\e_0[$ is chosen as above.
\smallskip\par\noindent
{\bf Claim:} For $x\in\Int(\DD_i)$, let $F^i_\e(x)$ denote the square obtained by the intersection of
${\vf_i}^{-1} ( \DD^i_\e )$ with the affine plane orthogonal to $\DD_i$ and containing $x$.
Then, for each $i=1,\ldots,m$ there exists a Borel set $\SS_i\sb\DD_i$, with $\sum_{i=1}^m\Ha^{n-2}(\SS_i)<c(n)\,\s$ for some absolute real constant
$c(n)>0$, such that for every $x\in \DD_i\sm\SS_i$, the 2-dimensional restriction of $u$ to the square $F^i_\e(x)$ is a Sobolev map with values into
$\Sph^1$ and with no homological singularities outside the center $x$ of the square.
The validity of this claim can be checked as a consequence of the mass estimate \eqref{sigma} and of a slicing argument.
\smallskip
\par We now modify the map $u$ as follows. For each $i=1,\ldots,m$, let $v_{i,\e}:\DD^i_\e\to\Sph^1$ be given by
\beq\label{ext}
v_{i,\e}(\wid x,\wih x):=u\Bigl( \e\,y_i(\wih x)\,\frac{\wid x}{\Vert\wid x\Vert} ,\wih x\Bigr)
\eeq
and define $u_\e:Q^n\to\Sph^1$ by
$$ u_\e(x)=\left\{ \ba{ll} v_{i,\e}({\vf_i}(x)) & {\textrm{if }} x\in {\vf_i}^{-1} ( \DD^i_\e )\,,\quad i=1,\ldots,m \\
u(x) & {\textrm{elsewhere in }}Q^n\,. \ea \right. $$
\par We can find a sequence $\{\e_h\}\searrow 0$ such that $u_h:=u_{\e_h}\in W^{1,1}(Q^n,\Sph^1)$ for each $h$, and $\{u_h\}$
converges to $u$ strongly in $W^{1,1}$.
%
The proof of this fact is omitted, since it follows by using arguments as in the next sections. To this purpose, we only observe that in Theorem~\ref{Tsimplex}, the computation in \eqref{conto gradiente} holds true even in case $k=2$.
Therefore, setting
$$ V_h:=Q^n\sm \bcup_{i=1}^m {\vf_i}^{-1} ( \DD^i_{\e_h} )\,,\quad h\in\Nat\,, $$
then $V_h$ is an open subset of $U$ with Lipschitz boundary (except at the points of the $(n-3)$-skeleton of $\rP$), and $\calL^n(V_h)\to\calL^n(U)=\calL^n(Q^n)$, as $h\to\i$.
\par Denote by $\rP_{u,h}$ the slice of the current $\rP(u)$ to the $(n-1)$-dimensional  boundary $\partial V_h$. Without loss of generality, we may and do choose the sequence $\{\e_h\}$ in such a way that $\rP_{u,h}$ is an $(n-3)$-rectifiable current
satisfying
\beq\label{Puh}
\e_h\,\gM(\rP_{u,h})\leq a_h\quad\fa\,h
\eeq
where $a_h\to 0$ as $h\to\i$.
\smallskip\par
We now apply the approximation theorem by Bethuel \cite[Thm.~2]{Bet}, see also \cite[Thm.~1.3]{HL}, to the Sobolev map
$u_{\vert V_{h}}:V_{h}\to\Sph^1$ where, we recall, $u_h=u$ on $V_h$.
This way, for each $h$ we find a sequence $\{v^{(h)}_k\}_k\sb R^\i(V_{h},\Sph^1)$ strongly converging to $u_{\vert V_{h}}$ in $W^{1,1}$.

Denote by $\rP(v^{(h)}_k)$ the $(n-2)$-current of the singularities of the Sobolev map $v^{(h)}_k$, so that $\spt \rP(v^{(h)}_k)\sb\bar V_h$ and
$(\pa \rP(v^{(h)}_k))\pri V_h=0$ for each $k$.
By an inspection to the construction of the approximating sequence from \cite{Bet,HL}, it turns out that
$$ \sup_k\gM(\rP(v^{(h)}_k))\leq c(n)\,\gM(\rP(u_h)\pri V_{h})\,. $$
In fact, the construction makes use of a slicing argument, where the degree of $v^{(h)}_k$ at the boundary of the 2-faces of the grid is bounded (up to an absolute constant) in terms of mass of the current of the singularities times $\d^{1-n}$, where $\d$ is the mesh of the grid. Therefore, the map $v^{(h)}_k$ being given by homogeneous extension on ``bad'' sets, one obtains the inequality in the last centered formula.

Therefore, since by \eqref{sigma} we can estimate
$$ \gM(\rP(u_h)\pri V_{\e_h})\leq\gM(\rP(u)\pri U)\leq \s\quad\fa\,h $$
we infer that
\beq\label{Pvh}
\sup_k\gM(\rP(v^{(h)}_k))\leq c(n)\,\s\quad\fa\,h\,.
\eeq
\par Moreover, viewing $\rP(v^{(h)}_k)$ as a current in $Q^n$, since the mass of the restriction of $\rP(v^{(h)}_k)$ to the boundary of $V_h$ is bounded (up to an absolute constant factor) in terms of the mass of the restriction of $\rP(u)$ to $\partial V_h$,
we also have:
$$ \sup_k\gM(\pa\rP(v^{(h)}_k))\leq c(n)\,\gM(\rP_{u,h})\,, $$
so that by \eqref{Puh} we can estimate
\beq\label{Phbis}
\sup_k\gM(\pa\rP(v^{(h)}_k))\leq c(n)\,\frac{a_h}{\e_h}\quad\fa\,h\,.
\eeq
\par In a way similar to definition \eqref{ext}, we now take for each $i$ the zero-homogeneous extension of $v^{(h)}_k$ in
${\vf_i}^{-1} ( \DD^i_{\e_h} )$ with respect to the coordinates $\wid x$ orthogonal to the $(n-2)$-simplex $\DD_i$.
We thus find a sequence $\{w^{(h)}_k\}\sb R^\i(Q^n,\Sph^1)$ strongly converging to $u$ in $W^{1,1}(Q^n,\gR^2)$  as $k$ tends to $\infty$.
Since
$$\gM(\rP(w^{(h)}_k))\pri (U\setminus V_h)\leq c(n)\,\e_h\,\gM(\pa\rP(v^{(h)}_k))\,, $$
by \eqref{Pvh} and \eqref{Phbis} we obtain the mass estimate
\beq\label{Pwhk}
\sup_k\gM\bigl((\rP(u)-\rP(w^{(h)}_k))\pri U\bigr)\leq c(n)\,(\s+a_h) \quad\fa\,h\,.
\eeq
%
%
\par By a diagonal argument, we thus find a sequence $\{w_h\}\sb R^\i(Q^n,\Sph^1)$ strongly converging to $u=u_\s$ in $W^{1,1}(Q^n,\gR^2)$ and such that
by \eqref{Pwhk} and \eqref{sigma}
\beq\label{Pwh}
\gM\bigl(\rP(w_h)\pri U\bigr)\leq c(n)\,(\s+a_h)
\quad\fa\,h\,.
\eeq
\par Now, recall that $\rP=\rP_\s$ satisfies \eqref{Pol}. By means of a slicing argument, we deduce that for $\Ha^{n-2}$-almost every $x\in{\textrm {int}}\,\DD_i$, where $i=1,\ldots,m$, the degree of $w_h$ around $x$ is a well-defined integer, that we denote by $d^{i}_h(x)$.
Therefore, we have:
\beq\label{Ph}
\rP(w_h)=\rP(w_h)\pri U+P_h
\eeq
where $P_h$ is an integral polyhedral chain with $\spt P_h\subset \spt\rP$, whose action is given by
$$ P_h(\y)=\sum_{i=1}^m
\int_{\DD_i}d^i_h\,\lan \y,\x_i\ran\,d\Ha^{n-2}\qquad\fa\,\y\in\D^{n-2}(Q^n) $$
where $\x_i$ is a unit $(n-2)$-vector orienting $\DD_i$, for $i=1,\ldots,m$.

Denote by $\theta$ the multiplicity of the current $\rP(u)$. By the previous Claim, we find a measurable set $K_\s$ contained in $\spt\rP$ such that $\Ha^{n-2}(K_\s)\leq c(n)\,\s$ and $$\sup_h|d^i_h(x)|\leq \theta(x) $$
for each $i=1,\ldots,m$ and $x\in \Int(\DD_i)\setminus K_\s$.
Moreover, we also estimate
$$ \sup_h\,\int_{K_\s}|d^i_h(x)|\,d\Ha^{n-2}(x)\leq c(n)\,E(\s)\,,\quad E(\s):=\gM(\rP(u_\s)\pri U) $$
so that $E(\s)\leq \s$.
Therefore, we obtain
$$ \sup_h\gM(P_h)\leq \gM(\rP_\s)+c(n)\,\s$$
and hence, on account of \eqref{Pwh} and \eqref{Ph},
$$ \gM(\rP(w_h))\leq c(n)\,(\s+a_h)+\gM(\rP_\s) \quad\fa\,h\,,$$
where, we recall, $a_h\to 0$ as $h\to\i$.
\par Letting $\s\searrow 0$ along a sequence, recalling that $\gM(\rP_\s)\to\gM(\rP(u))$,
by a further diagonal argument we find a sequence $\{u_k\}\sb W^{1,1}(Q^n,\Sph^1)$ strongly converging to $u$ in $W^{1,1}$ and such that
$\limsup_k\gM(\rP(u_k))\leq \gM(\rP(u))$. Since by lower semicontinuity $\gM(\rP(u))\leq \liminf_k\gM(\rP(u_k))$, we have proved Theorem~\ref{TBet}
under the assumption that $u$ is smooth near the boundary of $Q^n$.
\smallskip\par It the general case, we make use of a slicing argument as follows.
Let $\Vert x\Vert:=\sup_{1\leq i\leq n}|x_i|$ and $Q^n_\l=\{x\in\gR^n\,:\,\Vert x\Vert <\l\}$, so that $Q^n=Q^n_1$. For a.e. $0<\l<1$ the restriction $\rP(u)\pri Q^n_\l$ satisfies $\gM(\pa(\rP(u)\pri Q^n_\l))<\i$. Therefore, the boundary rectifiability theorem (cf. \cite[30.3]{Si}) implies that $\rP(u)\pri Q^n_\l$ is an integral $(n-2)$-current in $Q^n$, with support a closed set contained in the open cube $Q^n$.
We then apply again Federer's strong polyhedral approximation theorem, obtaining for each $\s>0$ small a
diffeomorphism $\vf_\s$ of $Q^n$ onto itself and an
$(n-2)$-dimensional polyhedral chain $\rP_\s$ with support in $Q^n$ such that $\vf_{\s\#}(\rP(u)\pri Q^n_\l)-\rP_\s=\pa R_\s$ for some current
$R_\s\in\R_{n-1}(Q^n)$ with $\gM(R_\s)+\gM(\pa R_\s)<\s$. Setting now $u_{\l,\s}(x):=u\circ{\vf_\s}^{-1}(x)$, $x\in \vf_\s(Q^n_\l)$, we have
$\rP(u_{\l,\s})=\vf_{\s\#}(\rP(u)\pri Q^n_\l)$. Therefore, arguing as before, we find a sequence $\{w_k\}\sb R^\infty(\vf_\s(Q^n_\l),\Sph^1)$ strongly converging to $u_{\l,\s}$ in $W^{1,1}(\vf_\s(Q^n_\l),\gR^2)$ and such that $\gM(\rP(w_k))\to\gM(\rP(u_{\l,\s}))$. Setting for $v=w_k$ or $v=u_{\l,\s}$
$$ \ol v(x):= v(\vf_\s(\l\,x))\,, \quad x\in Q^n \,, $$
and taking $\l\nearrow 1$, the claim follows through a diagonal argument.
\epf
\subsection{Energy approximation at the singular set}
\bpf[Proof of Theorem \ref{Tdipolo}.] Due to the condition on the degree around $\Delta$, setting for simplicity
\begin{align}\label{r eps}
r_\e(\wih x):= \e y(\wih x),
\end{align}
 we can find a smooth homotopy map $H:[0,1]\tim[0,2\pi]\tim\wih\DD\to\Sph^1$ such that
$H(0,\t,\wih x)=(\cos(d\theta),\,\sin(d\theta))$ and $H(1,\t,\wih x)=u(r_\e(\wih x)(\cos\t,\sin\t),\wih x)$. Since $u_{|\pa \DD_\e}\in W^{1,1}$, we can assume that $H\in W^{1,1}([0,1]\tim[0,2\pi]\tim\wih\DD,\Sph^1)$.
Define now in cylindrical coordinates $\wid x=\r(\cos\t,\sin\t)$ the map $v_\e:\DD_\e\to\gR^2$ as
$$\bar v_\e(\r,\t,\wih x):=\left\{ \ba{ll} H\Bigl( 2\r/r_\e(\wih x)-1,\t,\wih x\Bigr) & {\text{if }}
r_\e(\wih x)/2\leq |\wid x|\leq r_\e(\wih x) \\
2\r/r_\e(\wih x)\,(\cos(d\t),\sin(d\t)) & {\text{if }}
 |\wid x|\leq r_\e(\wih x)/2 \,.
\ea \right. $$
Notice that $v_\e$ is smooth on $\pa\DD_\e$ and $v_\e=u$ on $\pa\DD_\e$.
We claim that \eqref{est-dipolo} holds.

In fact, for $r_\e(\wih x)/2<\r< r_\e(\wih x)$, setting $t(\r)=2\r/r_\e(\wih x)-1$, we compute
$$ \pa_\r \bar v_\e(\r,\t,\wih x)=\pa_t H(t(\r),\t,\wih x)\cdot\frac 2{r_\e(\wih x)}\,,\qquad \pa_\t\bar v_\e(\r,\t,\wih x)=\pa_\t H(t(\r),\t,\wih x) $$
whereas
$$ \nabla_{\wih x} \bar v_\e(\r,\t,\wih x)=\nabla_{\wih x}H(t(\r),\t,\wih x)-\frac{2\r}{r_\e(\wih x)^2}\partial_t H(t(\r),\t,\wih x)\otimes\nabla r_\e(\wih x)\,. $$
We have 
$$
\begin{array}{l}
\ds \int\lm_{\DD_\e\cap\{r_\e(\wih x)/2<\r< r_\e(\wih x)\}}|\nabla v_\e|\,d\calL^n \\
\ds = \int_{\wih \DD}\int_0^{2\pi}\!\!\int_{r_\e(\wih x)/2}^{r_\e(\wih x)}\rho\sqrt{|\pa_\rho\bar v_\e|^2+\frac{|\pa_\theta\bar v_\e|^2}{\rho^2}+|\nabla_{\wih x}\bar v_\e|^2}d\rho d\theta d\wih x \\
\ds \leq\, \int_{\wih \DD}\int_0^{2\pi}\!\!\int_{r_\e(\wih x)/2}^{r_\e(\wih x)}\rho\Bigl[\frac{2}{r_\e(\wih x)}|\pa_tH|+\frac{|\pa_\theta H|}{\rho}+|\nabla_{\wih x}H|\\
\ds \qquad\qquad  +\frac{2\rho}{r_\e(\wih x)^2}|\pa_tH||\nabla r_\e| +\frac{2\sqrt \rho}{r_\e(\wih x)}\sqrt{|\nabla_{\wih x}H||\pa_tH||\nabla r_\e|}\Bigr]d\rho d\theta d\wih x,
\end{array}
$$
where all the partial derivatives of $H$ are computed at $(t(\r),\theta,\wih x)$ and $\nabla r_\e$ is computed at $\wih x$.
Using that $\rho\leq r_\e(\wih x)$ on $\DD_\e$ and $|\nabla r_\e(\wih x)|=\e$, for some absolute real constant $C$, we get
\begin{align*}
\int\lm_{\DD_\e\cap\{r_\e(\wih x)/2<\r< r_\e(\wih x)\}}|\nabla v_\e|\,d\calL^n
&\leq C\int_{\wih \DD}\int_0^{2\pi}\!\!\int_{r_\e(\wih x)/2}^{r_\e(\wih x)}|\nabla H(t(\r),\t,\wih x)|d\rho d\theta d\wih x\\
&= C\e\int\lm_{[0,1]\tim[0,2\pi]\tim\wih\DD}|\nabla H(t,\t,\wih x)|dt d\theta d\wih x=O(\e),
\end{align*}
where we performed the change of variable $t=t(\r)$ and we used that $r_\e(\wih x)\leq\e$.

On the other hand, by the area formula
$$\int\lm_{\DD_\e\cap\{r_\e(\wih x)/2<\r< r_\e(\wih x)\}}|M_2(\nabla v_\e)|\,d\calL^n=0\,. $$
Moreover, for $\r\leq r_\e(\wih x)/2$ we get
\begin{align*}
&\pa_\r\bar v_\e(\r,\t,\wih x)=\frac{2}{r_\e(\wih x)}(\cos(d\t),\sin(d\t)),\\
&\pa_\t\bar v_\e(\r,\t,\wih x)=\frac{2d\r}{r_\e(\wih x)}(-\sin(d\t),\cos(d\t)),\\
&\nabla_{\wih x}\bar v_\e(\r,\t,\wih x)=-\frac{4\r}{r_\e(\wih x)^2}(\cos(d\t),\sin(d\t))\otimes\nabla r_\e(\wih x).
\end{align*}
Therefore \begin{align*}
 \int\lm_{\DD_\e\cap\{\r<r_\e(\wih x)/2\}} \!\!\!\!\!\!
|\nabla v_\e|\,d\calL^n\leq&\int_{\wih \DD}\int_0^{2\pi}\!\!\int_{0}^{r_\e(\wih x)/2}\rho\left[\frac{2}{r_\e(\wih x)}+\frac{2|d|}{r_\e(\wih x)}+\frac{4\r\e}{r_\e(\wih x)^2}\right]d\r d\t d\wih x\\
&\leq2\pi\int_{\wih \DD}\int_0^{r_\e(\wih x)/2}\left[1+|d|+\frac{\e}{r_\e(\wih x)}\right]=O(\e),
\end{align*}
where we used that $|\nabla {r_\e(\wih x)}|=\e$ and $\r\leq{r_\e(\wih x)/2}$.
Finally we get \begin{align*}
\int\lm_{\DD_\e\cap\{\r<r_\e(\wih x)/2\}} \!\!\!\!\!\!
|M_2(\nabla v_\e)|\,d\calL^n&=\int_{\wih \DD}\int_0^{2\pi}\int_0^{r_\e(\wih x)/2}\r|M_2(\nabla \bar v_\e)|d\r d\t d\wih x\\
&\leq\int_{\wih \DD}\int_0^{2\pi}\int_0^{r_\e(\wih x)/2}\left[\frac{4|d|\r}{r_\e(\wih x)^2}+\frac{8|d|\r^2\e}{r_\e(\wih x)^2}\right]d\r d\t d\wih x\\
&=\int_{\wih \DD}\int_0^{2\pi}\frac{|d|}{2} d\t d\wih x + O(\e)\\
&=\pi\,|d|\,\Ha^{n-2}(\DD)+O(\e),
\end{align*}
so that \eqref{est-dipolo} readily follows. \epf
\subsection{Removal of point singularities}
\bpf[Proof of Theorem \ref{Tsing}.]
The argument being local, we may and do assume $\SS=\{0_{\gR^n}\}$. For $r>0$ small, we choose $v:B^n_r\to\gR^2$ smooth and such that $v=u$ on $\pa B^n_r$. We then define $w:Q^n\to\gR^2$ by taking
$$w(x)=\left\{ \ba{ll} u(x) & {\text{if }}|x|\geq r \\
u(rx/|x|) & {\text{if }}\d<|x|<r \\ v(rx/\d) & {\text{if }} |x|\leq\d \ea\right. $$
where $\d\in(0,r)$ is small, and we first estimate the energy of $w$ on $B^n_r\setminus B^n_\d$.

Denoting by $\nabla_\tau$ the tangential component of the derivative at $x\in\pa B^n_\r$, we get
$$
\left|\nabla w\left({x}\right)\right|=\frac r\rho\,\left|\nabla_\tau u\left(r\,\frac{x}{|x|}\right)\right|\,,\quad x\in B^n_r\sm B^n_\d  $$
and also
$$
\left|M_2(\nabla w\left({x}\right)\right)|=\Bigl(\frac r\rho\Bigr)^2\,\left|M_2\left(\nabla_\tau u\left(r\,\frac{x}{|x|}\right)\right)\right|\,,\quad x\in B^n_r\sm B^n_\d  $$
so that
$$
\begin{aligned}
\gM (G_w\pri (B^n_r\setminus B^n_\d)\tim\gR^2) = & \int_{B^n_r\setminus B^n_\d}\sqrt{1+|\nabla w|^2+|M_2(\nabla w)|^2}dx\\
\leq & |B^n_r|+\int_{B^n_r}|\nabla w|\,dx+\int_{B^n_r}|M_2(\nabla w)|\,dx \\
= &|B^n_r|+\int_0^r\frac r\rho\int_{\pa B^n_\rho}\left|\nabla_\tau u\left(r\,\frac{x}{|x|}\right)\right|\,d\Ha^{n-1}d\rho\\
& \ds
+\int_0^r\!\!\Bigl(\frac r\rho\Bigr)^2\!\!\int_{\pa B^n_\rho}\left|M_2\Bigl(\nabla_\tau u\Bigl(r\,\frac{x}{|x|}\Bigr)\Bigr)\right|\,d\Ha^{n-1}d\rho\\
=&|B^n_r|+\int_0^r\Bigl(\frac \rho r\Bigr)^{n-2}d\rho\int_{\pa B^n_r}\left|\nabla_\tau u(y)\right|\,d\Ha^{n-1} \\
& \ds +
\int_0^r\Bigl(\frac \rho r\Bigr)^{n-3}d\rho\int_{\pa B^n_r}\left|M_2(\nabla_\tau u(y))\right|\,d\Ha^{n-1}
\\
= &|B^n_r|+\frac r{n-1}\int_{\partial B^n_r}\left|\nabla_\tau u(y)\right|d\mathcal{H}^{n-1} \\
&\ds +\frac r{n-2}\int_{\partial B^n_r}\left|M_2(\nabla_\tau u(y))\right|d\mathcal{H}^{n-1}\,.
\end{aligned}
$$
Now, setting
$$F_1(r):=\int_{\partial B^n_r}\left|\nabla_\tau u\right|d\mathcal{H}^{n-1}, \quad F_2(r):=\int_{\partial B^n_r}\left|M_2(\nabla_\tau u)\right|d\mathcal{H}^{n-1}$$
we have
\begin{align*}
\int_0^1F_1(r)\,dr\leq\int_{B^n}|\nabla u|dx<\infty\,,\quad \int_0^1F_2(r)\,dr\leq\int_{B^n}|M_2(\nabla u)|dx<\infty\,.
\end{align*}
Thus we get necessarily that
$\liminf\lm_{r\to 0}r\,\bigl(F_1(r)+F_2(r)\bigr)=0$
and hence, definitely, 
$$
\liminf_{r\to0}\gM (G_w\pri (B^n_r\setminus B^n_\d)\tim\gR^2)=0\,.
$$
It remains to estimate the energy of $w$ on $B^n_\d$. We have
$$
\nabla w(x)=\frac{r}{\d}\nabla v\left(\frac{r}{\d}x\right),\qquad M_2(\nabla w)=\frac{r^2}{\d^2}M_2\left( \nabla v\left(\frac{r}{\d}x\right)\right)\,,\quad x\in B^n_\d\,.
$$
Then
$$
\begin{array}{l}
\qquad\qquad \ds\gM(G_w\pri B^n_\d\tim\gR^2) \leq \\
\ds  |B^n_\d|+\frac{r}{\d}\int_{B^n_\d}\left|\nabla v\left(\frac{r}{\d}x\right)\right|dx+\frac{r^2}{\d^2}\int_ {B^n_\d}\left|M_2\left( \nabla v\left(\frac{r}{\d}x\right)\right)\right|dx =\\
\ds  |B^n_\d|+\Bigl(\frac{\d}{r}\Bigr)^{n-1}\int_{B^n_r}|\nabla v(y)|dy+\Bigl(\frac{\d}{r}\Bigr)^{n-2}\int_ {B^n_r}|M_2( \nabla v(y))|dy<\infty\,.
\end{array}
$$
Therefore, recalling that $n\geq 3$, letting $r_j\to 0$ along a suitable sequence, and choosing $\d=\d(r_j)$ small w.r.t. $r_j$ we find $w_j:B^n_{r_j}\to\gR^2$ smooth with
$w_j=u$ on $\pa B^n_{r_j}$ such that:
$$ \lim_{j\to \infty}\gM(G_{w_j}\pri B^n_{r_j}\tim\gR^2)=0$$
and the proof is complete.
\epf
%
%
\subsection{Removal of high codimension singularities }
\bpf[Proof of Theorem \ref{Tsimplex}.]
Without loss of generality, we can assume that $\Sigma=\DD$, with $\DD$ an $(n-k)$-simplex, and that $\DD=\{0_{\gR^k}\}\times\wih \DD$, where we use the notation $x=(\wid x,\wih x)\in\gR^k\times\gR^{n-k}$.

Consider the neighborhood $\DD_\e=\{(\wid x,\wih x)\in \gR^n: \wih x\in\wih \DD\,,|\wid x|\leq r_\e(\wih x)\}$, where $r_\e$ is given by \eqref{r eps}, and for $r>0$ denote $\wid B_r:=\{\wid x\in \gR^k: |\wid x|\leq r\}$.

Let $\delta=\d(\e)<\e$ and $v:\DD_\e \to\gR^2$ be smooth such that
$v_{|\pa\DD_\delta}=u_{|\pa\DD_\delta}$, and define the map $w_\e:B^n\to\gR^2$ as
\begin{align*}
w_\e(x)=
\begin{cases}
u(x) &\mbox{in }B^n\setminus\DD_\e,\\
u\left(r_\e(\wih x){\wid x}/{|\wid x|},\wih x\right)\quad&\mbox{in }\DD_\e\setminus\DD_\delta,\\
v(\e \wid x/\delta,\wih x)&\mbox{in }\DD_\delta.
\end{cases}
\end{align*}
\par In order to estimate the energy on $\DD_\e\setminus\DD_\d$, let us start proving that
\begin{align}\label{grad w intercap}
\liminf_{\e\to0}\int_{\DD_\e\setminus\DD_\delta}|\nabla w_\e(x)|dx=0.
\end{align}
For this purpose, we introduce an adapted orthonormal basis $(\n,\tt_1,\ldots,\tt_{k-1})$ in $\gR^k$ so that $\n$ is the outward unit normal to $\pa \wid B_r$ at a point $\wid y\in\partial B_r$. Then, for any $x\in \DD_\e\setminus\DD_\d$, with
\beq\label{change}
\wid y=\wid y(x)=r_\e(\wih x)\,\frac{\wid x}{\r}\in\partial B_{r_\e(\wid x)}\,,\quad \r:=|\wid x|
\eeq
we obtain $ \pa_\n w_\e(\wid x,\wih x)=0$,
\beq\label{tau} \pa_{\tt_a}w_\e(\wid x,\wih x)=\frac{r_\e(\wih x)}\r\,\pa_{\tt_a}u(\wid y(x),\wih x)\,,\quad \a=1,\ldots,k-2
\eeq
and denoting $\wih x=(x_{k+1},\ldots,x_{n})$
\beq\label{ics}
\pa_{x_\be}w_\e(\wid x,\wih x)=\pa_{x_\be}u(\wid y(x),\wih x)+\pa_\n u(\wid y(x),\wih x)\,\pa_{x_\be}r_\e(\wih x)\,,\quad\be=k+1,\ldots,n\,.
\eeq
Therefore, since the distance function is 1-Lipschitz and $r_\e(\wih x)\leq\e$, there exists a positive constant $c$, only depending on $k$ and $n$, such that
$$ |\nabla w_\e(\wid x,\wih x)|\leq c\,\frac\e\r\,|\nabla u(\wid y(x),\wih x)|\,. $$
Using the change of variable in \eqref{change} and Fubini's theorem, we estimate:
\beq\label{conto gradiente}
\begin{array}{rl}
\ds \int\lm_{\DD_\e\setminus\DD_\delta}|\nabla w_\e(x)|\,dx \leq & \ds
c\,\e\int\lm_{\DD_\e}\r^{-1}\,|\nabla u(\wid y(x),\wih x)|\,dx \\
= & \ds c\,\e\,\int\lm_{\wih \DD}\Bigl(\int\lm_0^{r_\e(\wih x)}\r^{-1}\Bigl( \int_{\pa\wid B_\r} |\nabla u(\wid y(x),\wih x)|\,d\Ha^{k-1}\Bigr)\,
d\r\Bigr)\, d\wih x\\
= & \ds c\,\e\,\int\lm_{\wih \DD}\Bigl(\int\lm_0^{r_\e(\wih x)}\!\!\frac{\r^{k-2}}{r_\e(\wih x)^{k-1}}\Bigl( \int\lm_{\pa\wid B_{r_\e(\wih x)}}  \!\!\!\!\!\!
 |\nabla u(\wid y,\wih x)|\,d\Ha^{k-1}\Bigr)\,
d\r\Bigr)\, d\wih x \\
= & \ds \frac c{k-1}\,\e \int\lm_{\wih \DD}\Bigl( \int\lm_{\pa\wid B_{r_\e(\wih x)}}  |\nabla u(\wid y,\wih x)|\,d\Ha^{k-1}\Bigr)\,d\wih x \,,
\end{array}
\eeq
where we used that $k>2$.
Setting
$$ F(\e):=\int_{\wih \DD}\Bigl( \int_{\pa\wid B_{r_\e(\wih x)}}  |\nabla u(\wid y,\wih x)|\,d\Ha^{k-1}\Bigr)\,d\wih x $$
since for each $\e_0>0$ small
$$
\int_0^{\e_0} F(\e)\,d\e=\int_{\DD_{\e_0}}|\nabla u(x)|\,dx<\infty\,,
$$
then, necessarily $\liminf_{\e\to0}\e\, F(\e)=0$ and we obtain \eqref{grad w intercap}. \\
\par We now show that
\begin{align}\label{minor w intercap}
\liminf_{\e\to0}\int_{\DD_\e\setminus\DD_\delta}|M_2(\nabla w_\e(x))|dx=0.
\end{align}
We use again the adapted frame. This time, with an obvious notation, recalling that $ \pa_\n w_\e(\wid x,\wih x)=0$, by \eqref{tau} we get
for $ 1\leq\a_1<\a_2\leq k-1$
$$ |M_2(\nabla w_\e(\wid x,\wih x))_{\tt_{\a_1}\tt_{\a_2}}|=\Bigl(\frac{r_\e(\wih x)}{\r}\Bigr)^2|M_2(\nabla u(\wid y(x),\wih x))_{\tt_{\a_1}\tt_{\a_2}}|\,,
$$
whereas by \eqref{ics}, for each $k+1\leq\be_1<\be_2\leq n$ we can estimate
$$\ba{rl}
|M_2(\nabla w_\e(\wid x,\wih x))_{x_{\be_1}x_{\be_2}}|\leq  & |M_2(\nabla u(\wid y(x),\wih x))_{x_{\be_1}x_{\be_2}}| \\
& +
|\pa_{x_{\be_1}}r_\e(\wih x)|\,|M_2(\nabla u(\wid y(x),\wih x))_{\n x_{\be_2}}| \\
& +
|\pa_{x_{\be_2}}r_\e(\wih x)|\,|M_2(\nabla u(\wid y(x),\wih x))_{\n x_{\be_1}}|\,.
\ea $$
Finally, for each $\a=1,\ldots, k-1$ and $\be=k+1,\ldots,n$ we have:
$$ \begin{array}{l}
|M_2(\nabla w_\e(\wid x,\wih x))_{\tt_{\a}x_{\be}}|\leq \\
\ds\qquad \frac{r_\e(\wih x)}{\r}\,\Bigl(|M_2(\nabla u(\wid y(x),\wih x))_{\tt_{\a}x_{\be}}|
+|\pa_{x_{\be}}r_\e(\wih x)|\,|M_2(\nabla u(\wid y(x),\wih x))_{\tt_{\a}\n}| \Bigr)\,.
\end{array}
$$
We thus definitely obtain the upper bound:
$$ |M_2(\nabla w_\e(\wid x,\wih x))|\leq c\,\Bigl(1+\e+\frac {r_\e(\wih x)}\r+\Bigl(\frac{{r_\e(\wih x)}}{\r}\Bigr)^2 \Bigr) \, |M_2(\nabla u(\wid y(x),\wih x))| $$
for some absolute constant $c$, possibly depending on $\DD$.
Therefore, as in \eqref{conto gradiente}, this time we estimate:
$$\begin{array}{l}
\ds \int\lm_{\DD_\e\setminus\DD_\delta} 
|M_2(\nabla w_\e(\wid x,\wih x))|\,dx \leq  \ds
c 
\int\lm_{\DD_\e} \Bigl(1+\e+\frac {r_\e(\wih x)}\r+\Bigl(\frac{{r_\e(\wih x)}}{\r}\Bigr)^2 \Bigr) |M_2(\nabla u(\wid y(x),\wih x))|
\,dx \\
=  \ds c\int\lm_{\wih \DD} 
\int\lm_0^{r_\e(\wih x)}\Bigl(1+\e+\frac {r_\e(\wih x)}\r+\Bigl(\frac{{r_\e(\wih x)}}{\r}\Bigr)^2 \Bigr) 
\int\lm_{\pa\wid B_\r} |M_2(\nabla u(\wid y(x),\wih x))|
\,d\Ha^{k-1} 
d\r 
\, d\wih x\\
 =  \ds c
\int\lm_{\wih \DD}
\int\lm_{\pa\wid B_{r_\e(\wih x)}} 
 |M_2(\nabla u(\wid y,\wih x))|
\,d\Ha^{k-1}
d\wih x \cdot
\int\lm_0^{r_\e(\wih x)} 
\Bigl(1+\e+\frac {r_\e(\wih x)}\r+\Bigl(\frac{{r_\e(\wih x)}}{\r}\Bigr)^2\Bigr)
\Bigl( \frac{\r}{r_\e(\wih x)}\Bigr)^{k-1}
 d\r
\\
 \leq  \ds  c\,C(k)\,\e\, G(\e)
\end{array}
$$
for some dimensional constant $C(k)$, where we used again that $k>2$, $r_\e(\wih x)\leq\e$, and we have denoted
$$ G(\e):= \int_{\wih \DD}\Bigl( \int_{\pa\wid B_{r_\e(\wih x)}} |M_2(\nabla u(\wid y,\wih x))|
\,d\Ha^{k-1}\Bigr)\,. $$
Since again for $\e_0>0$ small
$$
\int_0^{\e_0} G(\e)\,d\e=\int_{\DD_{\e_0}}|M_2(\nabla u(x))|dx<\infty\,,
$$
then, necessarily $\liminf_{\e\to0}\e\, G(\e)=0$ and we obtain \eqref{minor w intercap}.
\par On the other hand, for almost every $x\in\DD_\delta$
$$
\nabla_{\wid x}w_\e(x)=\frac{\e}{\d}\nabla_{\wid x}v\left(\frac{\e}{\d}\wid x,\wih x\right),\qquad \nabla_{\wih x}w_\e(x)=\nabla_{\wih x}v\left(\frac{\e}{\d}\wid x,\wih x\right).
$$
whereas for $1\leq i<j\leq n$, with an obvious notation,
\begin{align*}
&M_2(\nabla w_\e)_{ij}=\frac{\e^2}{\d^2}M_2\left(\nabla v\left(\frac{\e}{\d}\wid x,\wih x\right)\right)_{ij}& &i,j\leq k,\\
&M_2(\nabla w_\e)_{ij}=\frac{\e}{\d}M_2\left(\nabla v\left(\frac{\e}{\d}\wid x,\wih x\right)\right)_{ij}& &i\leq k, \quad j> k,\\
&M_2(\nabla w_\e)_{ij}=M_2\left(\nabla v\left(\frac{\e}{\d}\wid x,\wih x\right)\right)_{ij}& &i,j>k,
\end{align*}
so that (recalling that $0<\d<\e$) definitely
$$ |\nabla w_\e(\wid x,\wih x)|\leq\frac{\e}{\d}\,|\nabla v(\wid y(x),\wih x)|\,,\quad |M_2(\nabla w_\e(\wid x,\wih x))|\leq\Bigl(\frac{\e}{\d}\Bigr)^2\,|M_2(\nabla v(\wid y(x),\wih x))|\,,  $$
where we have denoted
\beq\label{change2}
\wid y=\wid y(x)=\frac{\e}{\d}\wid x\,.
\eeq
Therefore, changing variable by \eqref{change2}, we get 
$$\begin{array}{l}
\gM(G_{w_\e}\pri\DD_\d\times\gR^2)\leq \\
\quad\leq  \ds |\DD_\d|+\frac{\e}{\d}\int_{\DD_\d}|\nabla v(\wid y(x),\wih x)|\,d\wid x\,d\wih x
+\Bigl(\frac{\e}{\d}\Bigr)^2\int_{\DD_\d}|M_2(\nabla v(\wid y(x),\wih x)|\,d\wid x\,d\wih x\\
\quad = \ds |\DD_\d|+\Bigl(\frac{\d}{\e}\Bigr)^{k-1}\int_{\DD_\e}|\nabla v(\wid y,\wih x)|\,d\wid x\,d\wih x
+\Bigl(\frac{\d}{\e}\Bigr)^{k-2}\int_{\DD_\e}|M_2(\nabla v(\wid y,\wih x)|\,d\wid x\,d\wih x \,.
\end{array}
$$
In conclusion, since $k\geq3$, letting $\e_j\to 0$ along a suitable sequence, and choosing $\d=\d(\e_j)$ small w.r.t. $\e_j$,
on account of \eqref{grad w intercap} and \eqref{minor w intercap}
we find
$$ \lim_{j\to \infty}\gM(G_{w_{\e_j}}\pri \DD_{\e_j}\tim\gR^2)=0$$
and the proof is complete. \epf
\section{Final remarks and open questions}\label{Sec:final}
In this final section, we briefly discuss whether our approach in terms of currents extends (with the expected modifications) to the wider class of maps $u\in\BV(B^n,\Sph^1)$ with finite relaxed energy \eqref{Area-rel}.
We then show that in case of both dimension and codimension at least equal to three, the corresponding relaxed area functional fails to be subadditive as a set function, even in the Sobolev case.
\subsection{The BV case.} The optimal lifting Cartesian current satisfies again \eqref{lift} (see \cite{Mu22}), where the distribution $\DDiv_{\bar\a}\gm_u$ is defined as in \eqref{Divmudef}, with an obvious extension of the adjoint notation to the $\gR^{2\tim n}$-valued measure $Du$.\\
In case $D^Ju=0$, recalling that $D_iu^j=\nabla_iu^j\,\calL^n+(D^Cu)^j_i$, we define the graph current $G_u$ in such a way that for every $(n-1)$-form $\o=\o^{(1)}$ as in \eqref{om}, this time we have
$$ G_u(\o^{(1)}) =  \ds\sum_{j=1}^2\sum_{i=1}^n
\ds\int\nolm_{B^n}\phi^j_i(x,u(x))\,d D_i
u^j\,. $$
It turns out that properties \eqref{bdGu} and \eqref{PuDivmu} continue to hold. Therefore, the lower bound given by Proposition~\ref{Plbd} readily extends. On the other hand, the upper bound inequality holds true provided that one is able to find a sequence
$\{u_k\}\sb W^{1,1}(B^n,\Sph^1)$ converging to $u$ strictly in the $\BV$-sense and satisfying $\lim_k\gM(\rP(u_k))=\gM(\rP(u))$, compare \cite{GM}.
In conclusion, we expect that our Main Result, Theorem~\ref{MainResult}, extends to the wider class of maps $u\in\BV(B^n,\Sph^1)$ such that $D^Ju=0$.
\smallskip
\par When $D^Ju\neq 0$, instead, even in low dimension $n=2$, we are very far from having an explicit formula, and even a characterization of the class of maps in $\BV(B^2,\Sph^1)$ with finite relaxed energy \eqref{Area-rel}.  The situation is much more complicate, since homological tools similar to the ones exploited in this paper fail to detect the energy gap.

For example, let $u\in \BV(B^2,\Sph^1)$ be the symmetric triple junction map, so that $u$ is constant in each third of the unit disk, with the three constant values $\alpha,\beta,\gamma$ equal to the vertices
of an equilateral triangle $T_{\alpha\beta\gamma}$ inscribed in the unit circle $\Sph^1$. According to \eqref{Divmudef} and using the decomposition formula \cite[(4.6)]{Mu22}, we get
$ |\DDiv\gm_u|(B^2)=|T_{\alpha\beta\gamma}|$, i.e. the area of the triangle $T_{\alpha\beta\gamma}$, recovering the exact energy gap (compare \cite{BSS}).
\par On the other hand, by slightly modifying Example 4.5 from \cite{BSS2}, so that the vertices of the two triangles in the target space belong to $\Sph^1$, one obtains a piecewise constant map in $\BV(B^2,\Sph^1)$,
with jump set equal to the union of twelve radii, in such a way that $|\DDiv\gm_u|(B^2)=0$,
but the energy gap is positive. By enforcing this modification also in Example 4.6 of \cite{BSS2}, one obtains even a piecewise constant map in $BV(B^2,\Sph^1)$ with infinite energy gap.  Therefore, the relevant (topological) singularity at the origin is not seen by any reasonable definition of the current of the singularities $\rP(u)$.
\subsection{A counterexample to the measure property}
We finally come back to the main conjecture in this framework: for each function $u\in \BV(B^n,\gR^2)$ with finite relaxed energy \eqref{Area-rel},
the localized functional $B\mapsto\ol \A_{\BV}(u,B)$ is subadditive on open sets, and hence it can be extended to a Borel measure on $B^n$.
That property should follow since the topology induced by strict convergence in $\BV$ is stronger enough, compared to the $L^1$-topology, see \eqref{L1_area}.
\par On the other hand, we now see that locality fails to hold when both dimension and codimension are strictly larger than two. This drawback is due to the fact that energy concentration in the relaxation process may occur on one-dimensional sets. Therefore, strict convergence fails to be strong enough in order to guarantee uniqueness of the Cartesian current enclosing the graph of Sobolev functions $u$ in $W^{1,1}(B^3,\gR^3)$. As in \cite{AcDM}, we build up our counterexample by means of the vortex map.
\par Denote by $\Sph^2$ the unit sphere in $\gR^3$, and let $u:B^3\to\Sph^2$ be given by $u(x)=x/|x|$. Then, $u\in W^{1,p}(B^3,\Sph^2)$ for each exponent $1\leq p<3$. Moreover, the cofactor function ${\mathrm{cof}\,}\nabla u$ belongs to $W^{1,q}(B^3,\gR^{3\tim 3})$ for $1\leq q<3/2$, and $\det\nabla u=0$ $\calL^3$-a.e. on $B^3$, by the area formula. Then,
for each open set $B\sb B^3$, the 3-dimensional area of the graph of the restriction $u_{\vert B}$ satisfies:
$$ \A(u,B)=\int_{B}\sqrt{1+|\nabla u|^2+|{\mathrm{cof}\,}\nabla u|^2}\,dx<\infty\,. $$
\par The graph 3-current $G_u$ is well-defined as before, in terms of the pull back of the graph map w.r.t. the approximate gradient, and actually $G_u$ is i.m. rectifiable in $\R_3(B^3\tim\gR^3)$, with finite mass $\gM(G_u)=\A(u,B^3)$, and a non zero boundary,
$$(\pa G_u)\pri B^3\tim\gR^3=-\d_{0_{\gR^3}}\times\qu{\Sph^2}\,, $$
see \cite[Sec. 3.2.2]{GMSl1}.
Roughly speaking, there are two qualitatively different ways
to fill the hole in the graph of $u$: inserting a ball
$\d_{0_{\gR^3}}\tim\qu{D^3}$, where $D^3$ is the (naturally oriented) unit ball in the target space, or a cylinder
$\qu{L}\tim\qu{\Sph^2}$, where $\qu{\Sph^2}:=\pa\qu{B^3}$ and $L$ is any oriented line
segment connecting a point in the boundary $\pa B^3$ of the
domain to the origin $0_{\gR^3}$. Therefore both the 3-currents $T_1$ and $T_2$,
$$ T_1:=G_u+\d_{0_{\gR^3}}\tim\qu{D^3}\,,\qquad T_2:=G_u+\qu{L}\tim\qu{\Sph^2}\,, $$
are Cartesian currents in $B^3\times\gR^3$.
Furthermore, it is not difficult to find two sequences $\{u_k^{(i)}\}\sb C^\infty(B^3,\gR^3)$, where $i=1,2$, satisfying the following properties:
\begin{itemize}
\item $u^{(i)}_k\to u$ strongly in $W^{1,1}(B^3,\gR^3)$, and hence strictly in $\BV$;
\item $G_{u^{(i)}_k}\wc T_i$ weakly in $\D_3(B^3\times\gR^3)$;
\item $\A(u^{(i)}_k)\to \gM(T_i)$, as $k\to \infty$, where
$$ \gM(T_1)=\A(u)+\frac{4\pi}3\,,\quad \gM(T_2)=\A(u)+4\pi\,. $$
\end{itemize}
\par
The smooth functions $u^{(1)}_k$ are equal to $x/|x|$ outside the ball $B^3_{1/k}$, where they cover once and with the appropriate orientation the ball $D^3$ in the target space.
\par
The smooth functions $u^{(2)}_k$, instead, take values in the unit sphere $\Sph^2$. They are equal to $x/|x|$ outside the ball $B^3_{1/k}$ and a small conical neighborhood $U_k$ of the segment $L$ with opening angle $1/k$, so that $\calL^3(L_k)\to 0$.
Moreover, they are smoothly extended on $U_k\setminus B^3_{1/k}$ in such a way that on each radius $r\in(1/k,1)$ they cover almost all the unit sphere $\Sph^2$. 
This can be done in such a way that the Brouwer degree of the smooth map ${u^{(2)}_k}_{\vert \pa B^3_r}:\pa B^3_r\to \Sph^2$ is equal to zero, for any $r\in(1/k,1)$. 
Therefore, each $u^{(2)}_k$ can be smoothly extended to the smaller ball $B^3_{1/k}$ by taking values in $\Sph^2$ and hence by paying a small amount of extra energy.
\br
In dimension $n=2$, instead of $n=3$, the analogous to sequence $\{u^{(2)}_k\}$ does not converge to the vortex map $u(x)=x/|x|$ in the strict $\BV$ sense, and hence in $W^{1,1}$, too. In fact, the gradient of $u^{(2)}_k$ in $U_k$ is of the order of $c_n k^{2-n}$ for some absolute constant $c_n>0$, and hence in case $n=2$ we get
$$ |Du^{(2)}_k|(B^2)=\int_{B^2}|\nabla u^{(2)}_k|\,dx \to \int_{B^2}|\nabla u|\,dx+c_2>|Du|(B^2)\,, $$
as $k\to\infty$. Whence, only $L^1$-convergence (or weak*-$\BV$ convergence) of $u^{(2)}_k$ to $u$ holds, see also \cite{Mu22}.
\er
\par By the previous construction, denoting for each open set $B\sb B^3$
$$ \ol\A_{\BV}(u,B):=\inf\Bigl\{\liminf_{k\to\i}\mathcal{A}(u_k,B)\mid \{u_k\}\sb C^1(B,\gR^3),\,\,
u_k\BVc u
\Bigr\}\,, $$
clearly $\ol\A_{\BV}(u,B^3)<\infty$. Let now $B^3_r$ be the open ball centered at the origin and with radius $r$. Using the sequence $\{u^{(1)}_k\}$ and a slicing argument, as e.g. in \cite[Lemma~5.2]{AcDM}, we can find a radius $r_3\in(0,1)$ such that if $r>r_3$, then
$$ \ol\A_{\BV}(u,B^3_r)=\A(u,B^3_r)+\frac{4\pi}3\,. $$
On the other hand, using the sequence $\{u^{(2)}_k\}$ we obtain for any $0<r\leq 1$ the inequality
$$ \ol\A_{\BV}(u,B^3_r)\leq \A(u,B^3_r)+4\pi\,r\,. $$
\par  In conclusion, arguing exactly as in \cite[Thm.~5.1]{AcDM}, it turns out that the localized functional $B\mapsto \ol\A_{\BV}(u,B)$ fails to be subadditive.
%


\subsection*{Acknowledgment}
We thank Andrea Marchese for useful discussions. The authors are members of the GNAMPA of INDAM.
%

\end{document}